\documentclass[10pt]{amsart}

\usepackage{amsmath, amssymb}
\usepackage{latexsym}
\usepackage{graphicx}
\usepackage{rotating}
\usepackage{comment}
\usepackage[all]{xy}
\usepackage{pgf}
\usepackage[T1]{fontenc}

\allowdisplaybreaks[1]

\numberwithin{equation}{section}

\begin{document}

\newcommand{\D}{\displaystyle}
\newcommand{\bo}{\bar{\omega}}
\newcommand{\bO}{\bar{\Omega}}
\newcommand{\bfo}{\boldsymbol{\omega}}
\newcommand{\bfO}{{\bf\Omega}}
\newcommand{\bfbo}{\bar{\boldsymbol{\omega}}}
\newcommand{\bfbO}{{\bf \bar{\Omega}}}
\newcommand{\bA}{{\bf A}}
\newcommand{\ba}{{\bf a}}
\newcommand{\bfS}{{\bf \mathcal{S}}}
\newcommand{\be}{\bar{\eta}}
\newcommand{\bfx}{{\bf x}}
\newcommand{\bfy}{{\bf y}}
\newcommand{\bfz}{{\bf z}}
\newcommand{\bfu}{{\bf u}}
\newcommand{\bfv}{{\bf v}}
\newcommand{\bfw}{{\bf w}}
\newcommand{\bfe}{{\bf e}}
\newcommand{\bff}{{\bf f}}
\newcommand{\bfg}{{\bf g}}
\newcommand{\bfG}{{\bf G}}
\newcommand{\bfgg}{\textbf{\textsc{g}}}
\newcommand{\bfGG}{{\bf \mathcal{G}}}
\newcommand{\bfeta}{{\bf \eta}}
\newcommand{\bfalpha}{{\bf \alpha}}
\newcommand{\tilo}{\tilde{\omega}}

\newcommand{\ms}{\negmedspace}

\newtheorem{theorem}{Theorem}
\newtheorem{corollary}{Corollary}
\newtheorem{definition}{Definition}
\newtheorem{example}{Example}

\title{Dynamic Equivalence of Control Systems via Infinite Prolongation}

\author{Matthew W. Stackpole}
\email{stackpol@colorado.edu}
\thanks{My thanks go out to Jeanne Clelland and George Wilkens.}

\subjclass[2000]{Primary }




\begin{abstract}
In this paper, we put issue of dynamic equivalence of control systems in the context of pullbacks of coframings on infinite jet bundles over the state manifolds. While much attention has been given to differentially flat systems, i.e. systems dynamically equivalent to linear control systems, the advantage of this approach is that it allowed us to consider control affine systems as well. Through this context we are able to classify all control affine systems of three states and two controls under dynamic equivalence of the type $(\bfx,\bfu)\mapsto \bfy(\bfx,\bfu)$.
\end{abstract}

\maketitle

\section{Introduction}

A control system is an underdetermined system of $n$ ordinary differential equations (ODEs),
\begin{eqnarray}
\label{control_system}
\dot{ {\bf x} } & = & {\bf f} ({\bf x} ,{\bf u} ).
\end{eqnarray}
Control systems show up in the design of electrical and mechanical systems, among other things. The variables ${\bf x}$ whose time evolution is determined by the ODEs are called state variables, while the ``free parameters'' ${\bf u}$ are called control variables. A control system can be viewed as a submanifold $\Sigma$ of the tangent bundle of the state space in the following way: given a manifold $M$ and a curve ${\bf x}:I\rightarrow M$, we say that ${\bf x}$ is a solution to the system $\Sigma\subset TM$ if $({\bf x}(t), \dot{\bf x}(t))$ lies in $\Sigma$ for all $t\in I$. The map $\mathbb{R}^s\rightarrow T_{\bf x}M$ given by ${\bf u}\mapsto \big( {\bf x}, {\bf f}({\bf x},{\bf u}) \big)$ is a parametrization of $\Sigma_{\bf x} = \Sigma\cap T_xM$ with the parameters ${\bf u}$ seen as local coordinates on $\Sigma_{\bf x}$.

A dynamic equivalence takes trajectories of one system, $\dot{\bf x} = {\bf f}({\bf x},{\bf u})$, to those of another, $\dot{\bf y} = {\bf g}({\bf y},{\bf v})$, and back again via maps between jet spaces which allow state derivatives to get mixed in:
\[
({\bf x},{\bf u},\dot{{\bf u}},\ldots,{\bf u}^{(J)}) \mapsto {\bf y}({\bf x},{\bf u},\dot{{\bf u}},\ldots,{\bf u}^{(J)}).
\]
Through the defining equation (\ref{control_system}), derivatives of state variables can be expressed in terms of control variables and their derivatives as well. Static (feedback) equivalence, which is a diffeomorphism of the state space, is a special case where ${\bf y}={\bf y}({\bf x})$.

Up to dynamic equivalence at the first jet level ($J=0$), i.e. ${\bf x}={\bf x}({\bf y},{\bf v})$ and ${\bf y}={\bf y}({\bf x},{\bf u})$, my results classify all control affine systems,
\[
\dot{{\bf x}} = {\bf f}^0({\bf x}) + u_i {\bf f}^i({\bf x}),
\]
of at most three states and two controls through the use of Cartan's method of equivalence. The main result of this paper is that every control affine system of three states and two controls falls into one of three classes under dynamic equivalence. The numbered rows represent these three classes. The entries in each row are systems that, while dynamically equivalent, are not statically equivalent.

\begin{center}
\begin{tabular}{|c|c|c|c|}
\hline
$1$ & $\dot{x}_1 = u_1$ & $\dot{x}_1 = u_1$ & $\dot{x}_1 = u_1$ \\
 & $\dot{x}_2 = u_2$ & $\dot{x}_2 = u_2$ & $\dot{x}_2 = u_2$ \\
 & $\dot{x}_3 = x_2$ & $\dot{x}_3 = x_2u_1$ & $\dot{x}_3 = 1 + x_2u_1$ \\
\hline
$2$ & $\dot{x}_1 = u_1$ & & \\
 & $\dot{x}_2 = u_2$ & & \\
 & $\dot{x}_3 = 0$ & & \\
\hline
$3$ & $\dot{x}_1 = u_1$ & & \\
 & $\dot{x}_2 = u_2$ & & \\
 & $\dot{x}_3 = 1$ & & \\
\hline
\end{tabular}
\end{center}
\section{Control Systems}
\label{controlsystems}

A system of ordinary differential equations (ODEs) with more variables than equations is called a \textbf{control system}. Locally a control system with $n+s$ variables and $n$ equations can be written in the form
\begin{eqnarray*}
\dot{x}_1 & = & f_1(x_1,\ \ldots,\ x_n,\ u_1,\ \ldots,\ u_s), \\
\dot{x}_2 & = & f_2(x_1,\ \ldots,\ x_n,\ u_1,\ \ldots,\ u_s), \\
 & \vdots & \\
\dot{x}_n & = & f_n(x_1,\ \ldots,\ x_n,\ u_1,\ \ldots,\ u_s). \\
\end{eqnarray*}
For our purposes, we will consider the functions $f_i$, $1\leq i\leq n$, to be $\mathcal{C}^\infty$.

Here, $x_i:\mathbb{R}\rightarrow \mathbb{R}$ and $u_j:\mathbb{R}\rightarrow \mathbb{R}$. We will use $t$ as our independent variable, and derivatives with respect to $t$ will be denoted by a dot: $\frac{dx_i}{dt} = \dot{x}_i$. This system of equations can be abbreviated with the single vector equation $\dot{\bfx} = \bff(\bfx,\bfu)$ where $\bfx = (x_1,\ \ldots,\ x_n)^T$, $\bfu = (u_1,\ \ldots,\ u_s)^T$, and $\bff = (f_1,\ \ldots,\ f_n)^T$. This type of control system is called \textbf{time independent} since there is no explicit $t$ dependence in the $f_i$.

In general, quantities that are vectors or matrices, like $\bfx$ above, will be written in bold face to distinguish them from scalars, like $x_i$.

The variables $x_i$ are known as the \textbf{state} variables, while the variables $u_j$ are known as the \textbf{control} variables. To explain the terminology, imagine a hovercraft on the surface of a two-dimensional lake. The state variables would be those needed to describe the state of the hovercraft on the lake: position of the hovercraft, which direction the hovercraft is turned, and the translational and rotational velocities of the hovercraft. The time evolution of state variables is predetermined, in this case by the Newton-Euler equations of motion, which are given explicitly in the example below. The control variables allow external influence of the state variables' time evolution. In the hovercraft scenario, control variables could describe the hovercraft's motor: the magnitude and direction of its thrust. Control variables are exactly what the hovercraft operator uses to control the system.
\section{Dynamic Equivalence}
\label{dynequiv}

Geometrically, a control system can be viewed as a submanifold $\overline{\Sigma}=\mathbb{R}\times\Sigma$ of $\mathbb{R}\times TM$ in the following manner: Given local coordinates $\bfx$ on $M$, the control system $\Sigma$ is a manifold with local coordinates $(\bfx,\bfu)$. With local coordinates $(\bfx,\dot{\bfx})$ on $TM$, there is an embedding
\begin{eqnarray*}
\iota:\mathbb{R}\times\Sigma\rightarrow\mathbb{R}\times TM
\end{eqnarray*}
given in coordinates by
\begin{eqnarray*}
(t,\bfx,\bfu)\mapsto (t,\bfx,\bff(\bfx,\bfu) ).
\end{eqnarray*}
This embedding $\iota$ pulls back the contact forms $\{ dx_i - \dot{x}_i\ dt\ |\ i=1,\ldots,n\  \}$ on $\mathbb{R}\times TM$ to the forms $\{\ \omega_i = dx_i - f_i(\bfx , \bfu )\ dt\ |\ i=1,\ldots,n\ \}$ on $\mathbb{R}\times\Sigma$.

A {\bf solution} to a control system also has a geometric interpretation. Let $\bfx(t)$ be a curve in $M$, i.e. $\bfx:\mathbb{R}\rightarrow M$, and define $p_1\bfx(t) = (\bfx(t),\dot{\bfx}(t))\in TM$. Such a curve $\bfx(t)$ is a solution to the control system $\Sigma$ if there exists a map ${\bf \sigma}:\mathbb{R}\rightarrow \Sigma$ that makes the following diagram commute:

\centerline{
\xymatrix{
\Sigma \ar@{->}[rr]^{\iota|_\Sigma} & & TM \ar@{->}[dd] \\
\\
\mathbb{R} \ar@{->}[rr]_{\bfx} \ar@{->}[uurr]^{p_1\bfx} \ar@{-->}[uu]^{\bf \sigma} & & M
}}
\noindent In particular,
\[
p_1\bfx(t) = \left( \iota|_\Sigma \circ {\bf \sigma} \right) (t),
\]
or in other words, $p_1\bfx(t)\in\Sigma$ for all $t$. Note that $\omega_i(\dot{{\bf \sigma}}) = 0$ for $i=1,\ldots,n$.

We will use the convention in this paper that a control system with a bar over it, for example $\overline{\Sigma}$, is a subbundle of $\mathbb{R}\times TM$, while a control system without the bar, $\Sigma$, is a subbundle of $TM$ which is the projection of $\overline{\Sigma}$. In fact, since we will be requiring that time be preserved through our equivalences, we will have $\overline{\Sigma} = \mathbb{R}\times\Sigma$.

\subsection{Jet Spaces}

Since the idea of dynamic equivalence is to allow a ``change of variables'' using higher order derivatives, we need a setting in which these higher order derivatives can be dealt with, much like the tangent bundle lets us work with first order derivatives. This setting is a jet space. We will say that curves $a,b:\mathbb{R}\rightarrow\mathbb{R}$ with $a(0)=b(0)=0$ have the same \textbf{$K$-jets} at $0$ if
\[
\frac{da}{dt}(0)=\frac{db}{dt}(0),\ \frac{d^2a}{dt^2}(0)=\frac{d^2b}{dt^2}(0),\ \ldots ,\ \frac{d^Ka}{dt^K}(0)=\frac{d^Kb}{dt^K}(0).
\]
Given $n$-dimensional differentiable manifolds $U$ and $V$ and maps $a,b:U\rightarrow V$ with $a(x)=b(x)=q$, we will say that $a$ and $b$ have the same $K$-jets at $x$ if for any differentiable maps $\phi:\mathbb{R}\rightarrow U$, $\psi:V\rightarrow\mathbb{R}$ with $\phi(0)=x$, $\psi\circ a\circ\phi$ and $\psi\circ b\circ\phi$ have the same $K$-jets at $0$.

Note that having the same $K$-jets at $x$ is an equivalence relation among maps from $U$ to $V$. Define the \textbf{$K^{th}$-order jet bundle of $M$}, denoted by $\mathcal{J}^K(M)$, to be the bundle over $M$ whose fiber $\mathcal{J}^K(M)_\bfx$ over a point $\bfx\in M$ is the space of curves $a:\mathbb{R}\rightarrow M$ modulo the equivalence relation of having the same $K$-jets at $\bfx$. Notice that with this definition, $\mathcal{J}^0(M) = M$ and $\mathcal{J}^1(M) = TM$, where the equality here is actually a bundle-preserving diffeomorphism.

Define the \textbf{prolongation map} $p_{j,k}$ which takes lifts of $\mathcal{C}^\infty$ curves from $M$ in $\mathcal{J}^j(M)$ to lifts of $\mathcal{C}^\infty$ curves from $M$ in $\mathcal{J}^{k}(M)$ ($j<k$) as follows.
\[
p_{j,k}(\bfx(t),\dot{\bfx}(t),\ddot{\bfx}(t),\ldots,\bfx^{(j)}(t)) = (\bfx(t),\dot{\bfx}(t),\ddot{\bfx}(t),\ldots,\bfx^{(j)}(t),\ldots,\bfx^{(k)}(t)\ )
\]
We will denote $p_{0,j}$ simply as $p_j$.

\subsection{Definition of Dynamic Equivalence}

Let $M$ and $N$ be smooth manifolds (state spaces) and 
\begin{eqnarray}
\label{cs_eq}
\begin{array}{ccc}
\Sigma\ :\ \dot{\bfx} & = & \bff(\bfx,\bfu) \\
\Lambda\ :\ \dot{\bfy} & = & \bfg(\bfy,\bfv)
\end{array}
\end{eqnarray}
control systems over their respective state spaces. 

We say control systems \eqref{cs_eq} on $M$ and $N$ are \textbf{dynamically equivalent} over open sets $\mathcal{U}\subset\mathcal{J}^{J+1}(M)$ and $\mathcal{V}\subset\mathcal{J}^{K+1}(N)$ for nonnegative integers $J$ and $K$ if there exist smooth maps $\Phi:\mathcal{U}\rightarrow N$ and $\Psi:\mathcal{V}\rightarrow M$ so that when restricted to the appropriate open sets:
\begin{enumerate}
\item for any solution $\bfx(t)$ of $\dot{\bfx}=\bff(\bfx,\bfu)$, $(\Phi\circ p_{J+1})(\bfx(t))$ is a solution to $\dot{\bfy} = \bfg(\bfy,\bfv)$,
\item for any solution $\bfy(t)$ of $\dot{\bfy}=\bfg(\bfy,\bfv)$, $(\Psi\circ p_{K+1})(\bfy(t))$ is a solution to $\dot{\bfx} = \bff(\bfx,\bfu)$,
\item the following diagram commutes for solutions,

\centerline{
\xymatrix{
& \mathcal{J}^{J+1}(M) \ar@{->}[dd] \ar@{->}[ddrr]_(0.25){\Phi} & & \mathcal{J}^{K+1}(N) \ar@{->}[dd] \ar@{-}[dl]^{\Psi} & \\
& & \ar@{->}[dl] & & \\
\mathbb{R} \ar@{->}[r]^\bfx \ar@{->}[uur]^{p_{J+1}\bfx} & \mathcal{J}^0(M) & & \mathcal{J}^0(N) & \mathbb{R} \ar@{->}[l]_{\bfy} \ar@{->}[uul]_{p_{K+1}\bfy}
}}
i.e. $\Psi\circ p_{K+1}\circ\Phi\circ p_{J+1}(x(t)) = x(t)$ whenever $x(t)$ is a solution of $\Sigma$, and $\Phi\circ p_{J+1}\circ\Psi\circ p_{K+1}(y(t)) = y(t)$ whenever $y(t)$ is a solution of $\Lambda$.
\end{enumerate}
Note that this means
\begin{eqnarray*}
\bfy & = & \Phi\left( \bfx,\dot{\bfx},\ldots,\bfx^{(J+1)} \right) , \\
\bfx & = & \Psi\left( \bfy,\dot{\bfy},\ldots,\bfy^{(K+1)} \right) .
\end{eqnarray*}

We will use the same notation for maps between jet spaces as we did for control systems, namely $\varphi: \mathcal{J}^{j}(M) \rightarrow \mathcal{J}^{k}(N)$ and $\bar\varphi: \mathbb{R}\times\mathcal{J}^{j}(M) \rightarrow \mathbb{R}\times\mathcal{J}^{k}(N)$ with $\bar\varphi = id \times\varphi$. Also note that in the definition of dynamic equivalence, we are using maps $\varphi: \mathcal{J}^{j}(M) \rightarrow \mathcal{J}^{k}(N)$, so they are defined in terms of the coordinates:
\[
(\ \bfx,\ \dot{\bfx},\ \ddot{\bfx},\ \ldots,\ \bfx^{(j)}\ ) \mapsto (\ \bfy,\ \dot{\bfy},\ \ddot{\bfy},\ \ldots,\ \bfy^{(k)}\ ).
\]
However, in practice we will be concerned only with the restrictions of these maps to the prolongations of control systems (defined below). Therefore, by way of the defining equations $\dot{\bfx} = \bff(\bfx,\bfu)$ and $\dot{\bfy} = \bfg(\bfy,\bfv)$ of the control systems, we will be looking at the restriction of $\varphi$ to the appropriate submanifolds with the following coordinates:
\[
(\ \bfx,\ \bfu,\ \dot{\bfu},\ \ldots,\ \bfu^{(j-1)}\ ) \mapsto (\ \bfy,\ \bfv,\ \dot{\bfv},\ \ldots,\ \bfv^{(k-1)}\ ).
\]

The proof of the following theorem should be clear from the definition, which is the same definition given in \cite{po}.

\begin{theorem}
Dynamic equivalence is an equivalence relation of control systems.
\end{theorem}

\textbf{Static (feedback) equivalence} is a special case of dynamic equivalence for which $J=K=-1$, i.e. $\Phi:M\rightarrow N$ is a diffeomorphism with $\Psi=\Phi^{-1}$. For static equivalence, we have $\bar\Phi_*\overline{\Sigma} = \overline{\Lambda}$ and $\bar\Psi_*\overline{\Lambda} = \overline{\Sigma}$. We say two systems are \textbf{locally static equivalent} over $\mathcal{U}\subset M$ and $\mathcal{V}\subset N$ if there exist coverings $\mathcal{U} = \bigcup_{\alpha\in A}\mathcal{U}^\alpha$ and $\mathcal{V} = \bigcup_{\alpha\in A}\mathcal{V}^\alpha$ such that the systems are static equivalent over each $\mathcal{U}^\alpha$ and $\mathcal{V}^\alpha$.

\begin{figure}[!h]
\begin{center}
\includegraphics[height=1.8in,width=2.5in]{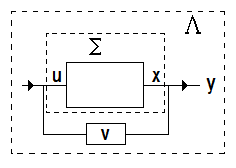}
\caption{Control system with feedback}
\label{system}
\end{center}
\end{figure}

From an engineering point of view, equivalence can be achieved through the addition of a feedback loop in the control system. In Figure \ref{system}, the system $\Sigma$ has input $\bfu$ and output $\bfx$. By adding a feedback loop, the new system $\Lambda$ has input $\bfv$ and output $\bfy$. In the case of static equivalence, the feedback loop only incorporates the old output $\bfx$ so that the new input $\bfv$ is only a function of $\bfx$ and $\bfu$. Including one or more integrators to the feedback loop allows $\bfv$ to be a function of $\bfx$, $\bfu$, and some number of derivatives of $\bfu$, and this is dynamic extension.

\subsection{Prolongation}

A key ingredient in dynamic equivalence is the notion of \textbf{prolongation of a control system}. For integers $k\geq 1$, define the \textbf{prolongation} of the system $\Sigma$ to the $k^{th}$ order, denoted by $\Sigma_k$, to be the subbundle of $\mathcal{J}^k(M)$ that corresponds to the prolongations of solutions of $\Sigma$, i.e. for any $x:I\rightarrow M$,
\begin{eqnarray*}
p_1(\bfx(t)) \in \Sigma\ \forall t\in I & \iff & p_k(\bfx(t)) \in \Sigma_k\ \forall t\in I.
\end{eqnarray*}
Obviously $\Sigma_1 = \Sigma$. In the same way that $\Sigma$ is a control system with $s$ control variables with state manifold $M$ of dimension $n$, we can view $\Sigma_2$ as a control system with $s$ control variables with state manifold $\Sigma$ of dimension $n+s$. An important fact is that $\Sigma$ is \textbf{strictly dynamically equivalent}, i.e. dynamically equivalent but not static equivalent, to $\Sigma_2$, as can be seen in the diagram below.

\centerline{
\xymatrix{
 & \Sigma^{2} \ar@{<-}[dd]_{p_{1,2}} \\
& \\
\Phi \ar@/_2pc/@{-}[r] & \Sigma \ar@{<-}[dd]_{p_{1}} \ar@/_2pc/@{<-}[l] \\
& \\
& \mathcal{J}^0(M) \ar@/_5pc/@{<-}[uuuu]_{\Psi}
}}

\begin{example}
The system
\begin{eqnarray*}
\Sigma:\ \dot{x}_1 & = & u_1 \\
\dot{x}_2 & = & u_2
\end{eqnarray*}
has two states and two controls. $\Sigma$ is dynamically equivalent to $\Sigma_2$:
\begin{eqnarray*}
\Sigma_2:\ \dot{y}_1 & = & y_3 \\
\dot{y}_2 & = & y_4 \\
\dot{y}_3 & = & v_1 \\
\dot{y}_4 & = & v_2
\end{eqnarray*}
where
\begin{eqnarray}
\label{spec_equiv}
\begin{array}{ccc}
x_1=y_1, &  x_2=y_2, & u_1=y_3, \\
u_2=y_4, & \dot{u}_1 = v_1, & \dot{u}_2 = v_2.
\end{array}
\end{eqnarray}
We have increased the number of states from two to four by viewing the controls as new state variables. \eqref{spec_equiv} gives the equivalence map. This is an example of what we will call a total prolongation.
\end{example}

In general, a \textbf{total prolongation} of the system $\dot{x} = f(x,u)$ is the system
\begin{eqnarray}
\label{prolongation}
\left( \begin{array}{c}
\dot{\bfx} \\
\dot{\bfu}
\end{array} \right)
 & = &
\left( \begin{array}{c}
\bff(\bfx,\bfu) \\
0
\end{array} \right)
+
\sum_i {\bf E}_i \dot{u}_i,
\end{eqnarray}
where ${\bf E}_i$ is the vector with a $1$ in the $(i+n)^{th}$ entry and zeros elsewhere. Here $(\bfx,\bfu)$ are the new state variables and $\dot{\bfu}$ are the new control variables. This system has a special form. A control system of the form
\begin{eqnarray}
\label{control_affine}
\dot{\bfx} = \bff(\bfx,\bfu) = \bff^0(\bfx) + \sum_i \bff^i(\bfx)u_i
\end{eqnarray}
is called \textbf{control affine}. In particular, \eqref{prolongation} is control affine. Thus we have the following theorem.

\begin{theorem}
\label{prolong}
Every control system $\Sigma$ is dynamically equivalent to a control affine system, namely $\Sigma_2$.
\end{theorem}

Similar to a total prolongation, some, but not all, of the control variables can be made into new state variables, as we see in the following example.

\begin{example}
The system
\begin{eqnarray*}
\Sigma:\ \dot{x}_1 & = & u_1 \\
\dot{x}_2 & = & u_2
\end{eqnarray*}
is dynamically equivalent to $\Lambda$:
\begin{eqnarray*}
\Lambda:\ \dot{y}_1 & = & y_3 \\
\dot{y}_2 & = & v_1 \\
\dot{y}_3 & = & v_2
\end{eqnarray*}
where
\begin{eqnarray*}
\begin{array}{ccc}
x_1=y_1, & x_2=y_2, & u_1=y_3, \\
 u_2=v_1, & \dot{u}_1 = v_2.
\end{array}
\end{eqnarray*}
We have increased the number of states from two to three by viewing only one of the controls as a new state variable. This process is called a \textbf{partial prolongation}. Every control system is dynamically equivalent to any partial prolongation of that system.
\end{example}

We will assume, without loss of generality, that the two systems in a dynamic equivalence have the same number of state variables ($m=n$). If $m<n$, perform repeated prolongations, either partial or total, until the number of states are equal and consider this new system.

A method for constructing a potential dynamic equivalence which is not a partial prolongation was mentioned briefly in a paper by Pomet \cite{po1}. Below we give a specific example of how the method works. This example incorporates both partial prolongation and changes of variables (a.k.a static equivalences) to give not only two control systems that are strictly dynamically equivalent but also the equivalence map.

\begin{example}
Start with an affine linear control system:
\begin{eqnarray}
\label{start}
\left( \begin{array}{c}
\dot{x}_1 \\
\dot{x}_2 \\
\dot{x}_3
\end{array}\right)
 & = &
\left( \begin{array}{c}
1 \\
0 \\
x_2
\end{array}\right)
u_1
+
\left( \begin{array}{c}
0 \\
1 \\
0
\end{array}\right)
u_2.
\end{eqnarray}
Partially prolong the three state system to a four state system.
\begin{eqnarray*}
z_1 = x_1 & \qquad & z_2 = x_2 \\
z_3 = x_3 & \qquad & z_4 = u_2 \\
w_1 = u_1 & \qquad & w_2 = \dot{u}_2
\end{eqnarray*}
\[
\left(\begin{array}{c}
\dot{z}_1 \\
\dot{z}_2 \\
\dot{z}_3 \\
\dot{z}_4 \\
\end{array}\right)
=
\left(\begin{array}{c}
0 \\
z_4 \\
0 \\
0 \\
\end{array}\right)
+
\left(\begin{array}{c}
1 \\
0 \\
z_2 \\
0 \\
\end{array}\right) w_1
+
\left(\begin{array}{c}
0 \\
0 \\
0 \\
1 \\
\end{array}\right) w_2
\]
By the nature of this partial prolongation, the $w_2$ vector must be of the form $(0\ 0\ 0\ 1)^T$. The systems $(\bfx,\bfu)$ and $(\bfz,\bfw)$ are dynamically equivalent. Through a change of basis, transform the $w_1$ vector into $(0\ 0\ 0\ 1)^T$.
\[
\left(\begin{array}{cccc}
-z_2 & -z_1 & 1 & 0 \\
0 & 1 & 0 & 0 \\
0 & 0 & 0 & 1 \\
1 & 0 & 0 & 0 \\
\end{array}\right)
\left(\begin{array}{c}
1 \\
0 \\
z_2 \\
0 \\
\end{array}\right)
=
\left(\begin{array}{c}
0 \\
0 \\
0 \\
1 \\
\end{array}\right)
\]
This corresponds to the change of coordinates
\[
(\tilde{z}_1,\tilde{z}_2,\tilde{z}_3,\tilde{z}_4) = (\ z_3-z_1z_2,\ z_2,\ z_4,\ z_1\ ).
\]
The change of coordinates is a static equivalence between $(\bfz,\bfw)$ and $(\tilde{\bfz},\bfw)$, and so we have yet another system dynamically equivalent to $(\bfx,\bfu)$.
\[
\left(\begin{array}{c}
\dot{\tilde{z}}_1 \\
\dot{\tilde{z}}_2 \\
\dot{\tilde{z}}_3 \\
\dot{\tilde{z}}_4 \\
\end{array}\right)
=
\left(\begin{array}{c}
-\tilde{z}_3\tilde{z}_4 \\
\tilde{z}_3 \\
0 \\
0 \\
\end{array}\right)
+
\left(\begin{array}{c}
0 \\
0 \\
0 \\
1 \\
\end{array}\right) w_1
+
\left(\begin{array}{c}
0 \\
0 \\
1 \\
0 \\
\end{array}\right) w_2
\]
The $(\tilde{\bfz},\bfw)$ will be a partial prolongation of a three state system. In this case,
\begin{eqnarray*}
\tilde{z}_1 = \tilde{y}_1, & \qquad & \tilde{z}_2 = \tilde{y}_3, \\
\tilde{z}_3 = \tilde{y}_2, & \qquad & \tilde{z}_4 = \tilde{v}_1, \\
w_1 = \dot{\tilde{v}}_1, & \qquad & w_2 = \tilde{v}_2.
\end{eqnarray*}
The numberings were chosen so that the final equations of the control system end up in this particularly nice form.
\begin{eqnarray*}
\left( \begin{array}{c}
\dot{\tilde{y}}_1 \\
\dot{\tilde{y}}_2 \\
\dot{\tilde{y}}_3
\end{array}\right)
& = & 
\left( \begin{array}{c}
0 \\
0 \\
\tilde{y}_2
\end{array}\right)
+
\left( \begin{array}{c}
-\tilde{y}_2 \\
0 \\
0
\end{array}\right) \tilde{v}_1
+
\left( \begin{array}{c}
0 \\
1 \\
0
\end{array}\right) \tilde{v}_2
\end{eqnarray*}

By construction, the systems $(\bfx,\bfu)$ and $(\tilde{\bfy},\bfv)$ are dynamically equivalent. What the process does not tell us is if this equivalence is strictly dynamic, for it could easily be static as well. In this example, however, the $(\bfx,\bfu)$ system is one of the classes of static equivalence given in Elkin \cite{el} and Wilkens \cite{wilk}, while the $(\tilde{\bfy},\bfv)$ system is clearly static equivalent to a distinct class $(\bfy,\bfv)$ 
\begin{eqnarray}
\label{end}
\begin{array}{ccc}
\left( \begin{array}{c}
\dot{y}_1 \\
\dot{y}_2 \\
\dot{y}_3
\end{array}\right)
& = & 
\left( \begin{array}{c}
0 \\
0 \\
y_2
\end{array}\right)
+
\left( \begin{array}{c}
1 \\
0 \\
0
\end{array}\right) v_1
+
\left( \begin{array}{c}
0 \\
1 \\
0
\end{array}\right) v_2
\end{array} 
\end{eqnarray}
following the transformation
\begin{eqnarray*}
\tilde{y}_i = y_i, & \quad & i=1,2,3, \\
v_1 = -\tilde{y}_2 \tilde{v}_1, & & \tilde{v}_2=v_2.
\end{eqnarray*}

It is interesting to note that unlike the original system \eqref{start} in our equivalence, \eqref{end} decouples into two smaller and separate systems: the first equation involves just $y_1$, $v_1$, while the other two equations involve only $y_2$, $y_3$, $v_2$. This equivalence also converts a nonlinear system $(\bfx,\bfu)$ into a linear one $(\bfy,\bfv)$. Both decoupling of equations and linearity greatly simplify the analysis of solutions of control systems.

Not only does the process presented above tell us that $(\bfx,\bfu)$ and $(\bfy,\bfv)$ are dynamically equivalent, but through some back tracking, it gives us the explicit equivalence maps.
\[
	(x_1,\ x_2,\ x_3,\ u_1,\ u_2,\ \ldots) \mapsto (x_3 - x_1x_2,\ u_2,\ x_2,\ -x_1u_2,\ \dot{u}_2,\ \ldots)
\]
\[
	(y_1,\ y_2,\ y_3,\ v_1,\ v_2,\ \ldots) \mapsto (-\frac{v_1}{y_2},\ y_3,\ y_1 - \frac{y_3v_1}{y_2},\ \frac{v_1v_2-y_2\dot{v}_1}{{y_2}^2},\ y_2,\ \ldots)
\]

This simple example also shows why it is necessary to consider dynamic equivalence on open sets. In this case, we would need to restrict our equivalence to the open set where $y_2\not=0$.

\end{example}
\section{Previous Results}
\label{previous}

The first theorem of this section is one of the most important, yet simplest to state, properties of dynamic equivalence. It can be found stated in a compatible form in \cite{fl}, but the following theorem and its proof, which are more in line with the terminology of this paper, can be found in \cite{po1}. 

\begin{theorem}
The number of control variables is an invariant of dynamic equivalence.
\end{theorem}

Note that while this theorem states that dynamically equivalent systems must have the same number of control variables, they may have different numbers of state variables. This is most obviously illustrated by Theorem \ref{prolong}. A system with $n$ states and $s$ controls is equivalent to its prolongation, which has $n+s$ states and $s$ controls. Thus the number of states in a system dynamically equivalent to a given system is unbounded.

Recall that a submanifold of an affine space is called \textbf{ruled} if, given any point of the submanifold, there is a line that passes through that point and that is contained completely within the submanifold. Classic examples of ruled submanifolds are planes, cylinders, and the hyperboloid of one sheet. We will abuse this terminology slightly and still call a submanifold ruled if it is the intersection of a ruled submanifold with a possibly bounded open set. A control system is called \textbf{ruled} if, when viewed as a subbundle $\Sigma$ of the tangent bundle $TM$, it defines at every point $x$ a ruled submanifold $\Sigma_x$ of the tangent space $T_xM$ at that point.

To state what is probably the most significant result to date in dynamic equivalence, some notation must be established. For $j<k$, let $\pi_{k,j}$ be the canonical projection from $\mathcal{J}^k(M)$ to $\mathcal{J}^j(M)$. Obviously $\pi_{k,k}$ is the identity. For any open set $\Omega\subset\mathcal{J}^k(M)$, define the subset $\Omega_l\subset\mathcal{J}^l(M)$ by
\begin{eqnarray*}
\Omega_l & = & \left\{\begin{array}{lc}
\pi_{k,l}(\Omega) & \textrm{if } l\leq k, \\
{\pi_{l,k}}^{-1}(\Omega) & \textrm{if } k\leq l.
\end{array} \right.
\end{eqnarray*}

The following is due to Pomet \cite{po}.

\begin{theorem} (Pomet)
\label{pomet_thm}
Let $\Sigma$ and $\Lambda$ be control systems with state manifolds $M$ and $N$ of dimension $m$ and $n$, $J$, $K$ two positive integers, and $\mathcal{U}\subset\mathcal{J}^{J+1}(M)$, $\mathcal{V}\subset\mathcal{J}^{K+1}(N)$ two open sets satisfying
\begin{eqnarray}
\label{open_requirement}
\mathcal{U}_1\cap\Sigma \subset (\mathcal{U}\cap\Sigma_{J+1})_1 & \textrm{ and } & \mathcal{V}_1\cap\Lambda \subset (\mathcal{V}\cap\Lambda_{K+1})_1.
\end{eqnarray}
If $\Sigma$ and $\Lambda$ are dynamic equivalent over $\mathcal{U}$ and $\mathcal{V}$, then 
	\begin{itemize}
	\item \textbf{if $m>n$}, then $\Sigma$ is ruled in $\mathcal{U}_1$.
	\item \textbf{if $n>m$}, then $\Lambda$ is ruled in $\mathcal{V}_1$.
	\item \textbf{if $m=n$}, then 
		\begin{itemize}
		\item (real analytic case) if $\mathcal{U}_1\cap\Sigma$ and $\mathcal{V}_1\cap\Lambda$ are connected, either $\Sigma$ and $\Lambda$ are ruled in $\mathcal{U}_1$ and $\mathcal{V}_1$, respectively, or they are locally static equivalent over $\mathcal{U}_1$ and $\mathcal{V}_1$.
		\item ($\mathcal{C}^\infty$ case) there are open sets $R,S\subset\mathcal{U}_1$ and $\mathcal{R},\mathcal{S}\subset\mathcal{V}_1$ with
			\begin{enumerate}
				\item $\mathcal{U}_1 = \bar{R}\cup S = R\cup\bar{S}$,
				\item $\mathcal{V}_1 = \bar{\mathcal{R}}\cup \mathcal{S} = \mathcal{R}\cup\bar{\mathcal{S}}$,
				\item $\Sigma$ and $\Lambda$ are ruled over $R$ and $\mathcal{R}$,
				\item $\Sigma$ and $\Lambda$ are static equivalent over $S$ and $\mathcal{S}$.
			\end{enumerate}
		\end{itemize}
	\end{itemize}
\end{theorem}

The condition \ref{open_requirement} basically says that nothing is lost when either prolonging the control system up or projecting the open set down in the jet spaces. In fact this containment is an equality; the reverse inclusion follows directly from the definitions.

Recall that every system $\Sigma$ is dynamically equivalent to its prolongation $\Sigma_2$. Since the dimension of the state space of $\Sigma_2$ is larger than the dimension of the state space of $\Sigma$, this theorem guarantees that $\Sigma_2$ is ruled. Of course we already know that $\Sigma_2$ is affine linear, so in this case the result is trivial. A natural question that arises from this, and one that is partially answered by this paper, is this:

\begin{center}
Given an affine linear control system, when is it the prolongation of a smaller system?
\end{center}

At the moment, this question has not been answered in its full generality, here or elsewhere. In an attempt to partially address this issue, this paper will classify control systems of low dimension that are affine linear up to dynamic equivalence in Chapter \ref{affinesystems}. The methods used to do this rely on a previous classification of affine linear control systems under static equivalence. Gardner and Shadwick first classified control systems with two state variables and one control variable under static equivalence \cite{feed1}. Wilkens then solved the problem for three states and two controls \cite{gw}. For the complete classification of affine linear systems under static equivalence with at most three states, which I present here without proof, see \cite{el}.

In the following theorem, $n$ represents the number of state variables $x_i$, and $u_j$ are control variables. Given a control system $\Sigma:\ \dot{\bfx} = \bff(\bfx,\bfu)$ with state space $M$, we say that a point $p\in M$ is \textbf{regular} if there is a neighborhood of $p$ on which the rank of $\Sigma$, defined to be the rank of $\frac{\partial\bff}{\partial\bfu}$, is constant.

\begin{theorem} (Elkin)
An affine linear control system \eqref{control_affine} with $n\leq 3$ states is locally static equivalent at a regular point $p$ to one of the following systems:
\begin{itemize}
	\item $n=1$
		\begin{eqnarray*}
		\dot{x}_1 = 0, & \dot{x}_1 = 1, & \dot{x}_1 = u_1.
		\end{eqnarray*}
	\item $n=2$
		\begin{eqnarray*}
			\left\{\begin{array}{ccc}
			\dot{x}_1 & = & 0 \\
			\dot{x}_2 & = & 0
			\end{array}\right. ,
			&
			\left\{\begin{array}{ccc}
			\dot{x}_1 & = & 1 \\
			\dot{x}_2 & = & 0
			\end{array}\right. ,
			&
			\left\{\begin{array}{ccc}
			\dot{x}_1 & = & u_1 \\
			\dot{x}_2 & = & 0
			\end{array}\right. ,
			\\
			\left\{\begin{array}{ccc}
			\dot{x}_1 & = & u_1 \\
			\dot{x}_2 & = & 1
			\end{array}\right. ,
			&
			\left\{\begin{array}{ccc}
			\dot{x}_1 & = & u_1 \\
			\dot{x}_2 & = & x_1
			\end{array}\right. ,
			&
			\left\{\begin{array}{ccc}
			\dot{x}_1 & = & u_1 \\
			\dot{x}_2 & = & u_2
			\end{array}\right. ,
		\end{eqnarray*}
	\item $n=3$
		\begin{eqnarray*} \begin{array}{c@{}c@{}c@{}c}
			\left\{\begin{array}{ccc}
			\dot{x}_1 & = & 0 \\
			\dot{x}_2 & = & 0 \\
			\dot{x}_3 & = & 0
			\end{array}\right. ,
			&
			\left\{\begin{array}{ccc}
			\dot{x}_1 & = & 1 \\
			\dot{x}_2 & = & 0 \\
			\dot{x}_3 & = & 0
			\end{array}\right. ,
			&
			\left\{\begin{array}{ccc}
			\dot{x}_1 & = & u_1 \\
			\dot{x}_2 & = & 0 \\
			\dot{x}_3 & = & 0
			\end{array}\right. ,
			&
			\left\{\begin{array}{ccc}
			\dot{x}_1 & = & u_1 \\
			\dot{x}_2 & = & 1 \\
			\dot{x}_3 & = & 0
			\end{array}\right. ,
			\\
			\left\{\begin{array}{ccc}
			\dot{x}_1 & = & u_1 \\
			\dot{x}_2 & = & x_1 \\
			\dot{x}_3 & = & 0
			\end{array}\right. ,
			&
			\left\{\begin{array}{ccc}
			\dot{x}_1 & = & u_1 \\
			\dot{x}_2 & = & x_1 \\
			\dot{x}_3 & = & 1
			\end{array}\right. ,
			&
			\left\{\begin{array}{ccc}
			\dot{x}_1 & = & u_1 \\
			\dot{x}_2 & = & x_1 \\
			\dot{x}_3 & = & x_2
			\end{array}\right. ,
			&
			\left\{\begin{array}{ccc}
			\dot{x}_1 & = & u_1 \\
			\dot{x}_2 & = & H(x)u_1 \\
			\dot{x}_3 & = & 1\! +\! x_2u_1
			\end{array}\right. ,
			\\
		\end{array} \end{eqnarray*}
		where $H(x)$ is an arbitrary function with $\frac{\partial H}{\partial x_3}$ is nonzero.
		\begin{eqnarray*}
			\left\{\begin{array}{ccc}
			\dot{x}_1 & = & u_1 \\
			\dot{x}_2 & = & u_2 \\
			\dot{x}_3 & = & 0
			\end{array}\right. ,
			&
			\left\{\begin{array}{ccc}
			\dot{x}_1 & = & u_1 \\
			\dot{x}_2 & = & u_2 \\
			\dot{x}_3 & = & 1
			\end{array}\right. ,
			&
			\left\{\begin{array}{ccc}
			\dot{x}_1 & = & u_1 \\
			\dot{x}_2 & = & u_2 \\
			\dot{x}_3 & = & u_3
			\end{array}\right. ,
			\\
			\left\{\begin{array}{ccc}
			\dot{x}_1 & = & u_1 \\
			\dot{x}_2 & = & u_2 \\
			\dot{x}_3 & = & x_2
			\end{array}\right. ,
			&
			\left\{\begin{array}{ccc}
			\dot{x}_1 & = & u_1 \\
			\dot{x}_2 & = & u_2 \\
			\dot{x}_3 & = & x_2u_1
			\end{array}\right. ,
			&
			\left\{\begin{array}{ccc}
			\dot{x}_1 & = & u_1 \\
			\dot{x}_2 & = & u_2 \\
			\dot{x}_3 & = & 1+x_2u_1
			\end{array}\right. .
		\end{eqnarray*}
\end{itemize}
\end{theorem}

The following theorem takes the problem of classifying control systems with one control variable under dynamic equivalence and reduces it to the simpler case of static equivalence. While this theorem has been known for some time (see \cite{po1} for one example), a new proof of this theorem in the framework of this paper will be given in Chapter \ref{scalarcontrol}.

\begin{theorem}
Let the control systems $\Sigma$, $\Lambda$ in \eqref{cs_eq} be dynamically equivalent with $s=1$ control variable and $m$, $n$ state variables, respectively. If $m=n$, then the systems are in fact static equivalent. If $m<n$ ($m>n$), then the systems are static equivalent after a finite number of prolongations of the smaller system $\Sigma$ ($\Lambda$).
\end{theorem}
\section{The Equivalence Problem}
\label{equivalence}

Given a manifold $M$, a \textbf{framing} on $M$ is a collection $\{ {\bf X}_i \}_{i=1}^n$ of smooth sections of the tangent bundle $TM$ such that for every $p\in M$, the collection of vectors $\{ ({\bf X}_i)_p \}_{i=1}^n$, called a \textbf{frame}, forms a basis for $T_pM$. A \textbf{coframing} is simply the dual of this notion, i.e. a collection of 1-forms $\{ \omega^j \}_{j=1}^n$ (smooth sections of the cotangent bundle $T^*M$) such that $\{ (\omega^j)_p \}_{j=1}^n$ forms a basis for $T^*_pM$ for every $p\in M$. Every coframing $\omega^j$ has a corresponding framing $X_i$ for which $\omega^j({\bf X}_i) = \delta^j_i$. A \textbf{local framing/coframing} is simply a framing/coframing defined on an open set $U\subset M$.

An \textbf{equivalence problem} \cite{moe} can be stated in the following way: Let $M^n$ and $N^n$ be smooth $n$-dimensional manifolds and $G\subset GL(n,\mathbb{R})$ a subgroup. Let $\omega_U = (\omega^1_U,\ldots,\omega^n_U)^T$ and $\Omega_V = (\Omega^1_V,\ldots,\Omega^n_V)^T$ be local coframings of $U\subset M$ and $V\subset N$, respectively, chosen in some geometrically natural way. We wish to find necessary and sufficient conditions that there exists a diffeomorphism $\varphi:U\rightarrow V$ such that
\begin{eqnarray*}
\varphi^*\Omega_V & = & \gamma_{VU}\omega_U
\end{eqnarray*}
where $\gamma_{VU}:U\rightarrow G$. A common abuse of notation, one which will be used in this paper, is to drop the pullback from the notation where the map $\varphi$ is clear from context: $\Omega_V = \gamma_{VU}\omega_U$.

For example, suppose we are given manifolds $M$ and $N$ with Riemannian metrics $ds^2$ and $dS^2$, respectively. We can locally diagonalize the metrics on open sets $U\subset M$ and $V\subset N$ such that
\begin{eqnarray*}
ds^2 = \sum_i (\omega^i_U)^2, & \qquad & dS^2 = \sum_i (\Omega^i_V)^2.
\end{eqnarray*}
The problem then is to find necessary and sufficient conditions such that a diffeomorphism $\varphi:M\rightarrow N$ exists such that $\varphi^*\Omega_V = \gamma_{VU}\omega_U$, where $\gamma_{VU}$ is an element of the orthogonal group $O(n)$.

The goal of this paper is to adapt the framework of an equivalence problem to dynamic equivalence. Then, using methods of exterior differential systems, we will classify a collection of control systems. What makes the dynamic equivalence problem tricky is the unboundedness of the size of the potentially equivalent state manifold, and hence also the lack of diffeomorphisms. A diffeomorphism $\varphi:M\rightarrow N$ cannot exist due to differences in dimension. In fact, strict dynamic equivalences are defined in terms of submersions rather than diffeomorphisms. This difficulty due to submersions persists through any finite number of prolongations. To solve this problem with submersions, in the next section we will simply make everything the same size: infinite.
\section{Infinite Prolongations}
\label{infinite}

The trick to dealing with our submersion woes is through prolongation, an idea introduced in section \ref{dynequiv}. Recall that a control system on $M$
\begin{eqnarray}
\label{f(x,u)}
\overline{\Sigma} : & \dot{x}_i=f_i(\bfx,\bfu), & \quad 1\leq i\leq n,
\end{eqnarray}
can be represented by
\[
{\bf X} = \frac{\partial}{\partial t} + \sum_{i=1}^n f_i(\bfx,\bfu)\frac{\partial}{\partial x_i}
\]
as a parametrization of $\overline\Sigma$ inside $\mathbb{R}\times TM$. A basis for the space ${\bf X}^\perp$ is
\begin{eqnarray*}
\omega_i & = & dx_i - f_i(\bfx,\bfu)\ dt, \quad 1\leq i\leq n.
\end{eqnarray*}
The forms $\omega_i$ are the pullback to $\overline{\Sigma}$ by the inclusion map of the contact forms $dx_i - \dot{x}_i\ dt$ on $\mathbb{R}\times TM$, where $\mathbb{R}\times TM$ has coordinates $(t,x_i,\dot{x}_i)$.The collection of $1$-forms $\{ dt,\omega^i,du_j \}$ forms a coframing on $\overline\Sigma$ that encodes the information of the control system.

Prolongation of \eqref{f(x,u)} yields a system $\overline\Sigma_2$ given by the equations
\begin{eqnarray*}
\dot{\bfx} & = & \bff(\bfx,\bfu) \\
\dot{\bfu} & = & \bar{\bfu}
\end{eqnarray*}
with state variables $\bfx,\bfu$ and control variables $\bar{\bfu}$. Thus a suitable coframing on $\overline{\Sigma}_2$ that encodes the information of the original control system and $\overline{\Sigma}_2$ is
\[
\left\{\begin{array}{ccc}
\omega^{-1} & = & dt, \\
\omega_i^0 & = & dx_i - f_i(\bfx,\bfu)\ dt, \quad 1\leq i\leq n, \\
\omega_j^1 & = & du_j - \bar{u}_j\ dt, \quad 1\leq j\leq s, \\
\omega_j^2 & = & d\bar{u}_j, \quad 1\leq j\leq s. \end{array}\right\}
\]

Define the \textbf{infinite jet bundle} $\mathcal{J}^\infty(M)$ as the projective limit of the finite jet bundles $\D \mathcal{J}^\infty(M) = \varprojlim_K \mathcal{J}^K(M)$, endowed with the projective limit topology. Let $\overline{\Sigma}_\infty$ and $\overline{\Lambda}_\infty$ be the projective limits of the prolongations of the control systems $\overline{\Sigma}$ and $\overline{\Lambda}$, respectively. By repeated iterations of the prolongation process above, a suitable choice for preferred coframings on $\overline{\Sigma}_\infty$ and $\overline{\Lambda}_\infty$ with coordinates $(t,\bfx,\bfu,\dot{\bfu},\ddot{\bfu},\ldots)$ and $(t,\bfy,\bfv,\dot{\bfv},\ddot{\bfv},\ldots)$, respectively, which encodes the information of the respective control systems is as follows.

\begin{eqnarray}
\label{coframe}
\begin{array}{c}
\bfo = 
\left( \begin{array}{c}
\omega^{-1} \\
\bfo^0 \\
\bfo^{1} \\
\bfo^{2} \\
\vdots
\end{array} \right)
=
\left( \begin{array}{c}
dt \\
d{\bf x} - f({\bf x},{\bf u})dt \\
d{\bf u} - \dot{{\bf u}} dt \\
d\dot{{\bf u}} - \ddot{{\bf u}} dt \\
\vdots
\end{array} \right)\\
\\
\bfO = 
\left( \begin{array}{c}
\Omega^{-1} \\
\bfO^0 \\
\bfO^{1} \\
\bfO^{2} \\
\vdots
\end{array} \right)
=
\left( \begin{array}{c}
dt \\
d{\bf y} - g({\bf y},{\bf v})dt \\
d{\bf v} - \dot{{\bf v}} dt \\
d\dot{{\bf v}} - \ddot{{\bf v}} dt \\
\vdots
\end{array} \right)
\end{array} 
\end{eqnarray}
The covectors $\bfo^i$ and $\bfO^i$ are $n$-dimensional for $i=0$ and $s$ dimensional for $i>0$.

Now we should take a closer look at what happens to the mappings involved in the definition of dynamic equivalence under this infinite prolongation process. Given a map $\Phi:\Sigma_{J+1}\rightarrow N$, as in the definition of dynamic equivalence, define the \textbf{$k^{th}$ prolongation} of the map, denoted $\Phi_{[k]}$ as the map that makes the following diagram commute on solutions. \\
\centerline{
\xymatrix{
& \Sigma_{J+k+1} \ar@{->}[dd] \ar@{->}[drr]^(0.5){\Phi_{[k]}} \\
& & & \Lambda_k \ar@{<-}[dd]_{p_k} \\
& \Sigma_{J+1} \ar@{->}[d] \ar@{->}[drr]^(0.5){\Phi} \\
I \ar@{->}[r]^{x} \ar@{->}[ur]^{p_{J+1}x} \ar@/^3pc/@{->}[uuur]^{p_{J+k+1}x} & M & & N
}} \\
In other words,
\begin{eqnarray*}
\left( p_k \circ \Phi \right)\big(\ p_{J+1} x(t) \ \big) & = & \Phi_{[k]} \big(\ p_{J+k+1} x(t) \ \big)
\end{eqnarray*}
for solutions $x(t)\in M$.

Now define $\Phi_\infty:\Sigma_\infty\rightarrow\Lambda_\infty$ by $\Phi_\infty = \varprojlim_{k}\Phi_{[k]}$ in the obvious fashion, i.e. for projection the projection map $\pi_k$ that takes an infinite jet to the $k^{\textrm{th}}$ jet,
\[
\pi_{k}\circ\Phi_\infty = \Phi_{[k]}\circ\pi_{J+k+1}.
\]
Let $\Psi:\Lambda_{K+1}\rightarrow M$ be the map used in section \ref{dynequiv} in the definition of dynamic equivalence, and define $\Psi_\infty$ similarly. From the definitions of dynamic equivalence and prolongation, it is simple to show that
\[
\Psi\circ\Phi_{[K+1]}\circ p_{J+K+2} = Id_0
\]
is the identity on curves in $M$. Finite prolongation of this relation shows
\begin{eqnarray}
\label{inverses}
\Psi_{[k]}\circ\Phi_{[K+1+k]}\circ p_{k,J+K+2+k} = Id_k
\end{eqnarray}
is the identity on curves in $\Sigma_k$. Taking the limit of \eqref{inverses} as $k$ tends to infinity tells us that
\[
\Psi_{\infty}\circ\Phi_{\infty} = Id_\infty
\]
is the identity on $\Sigma_\infty$. Similarly
\[
\Phi_{\infty}\circ\Psi_{\infty} = Id_\infty
\]
is the identity on $\mathcal{J}^\infty(N)$, and we can conclude that ${\Phi_\infty}^{-1} = \Psi_\infty$ and that $\Phi_\infty$ is a diffeomorphism.

To recap, in order to pose an equivalence problem for dynamic equivalence, we need a diffeomorphism between spaces. The problem with dynamic equivalence is that the maps used in the definition of the equivalence can never give us a diffeomorphism at any finite level (unless the equivalence is actually static). By passing to the infinite prolongation, the submersions become diffeomorphisms. We obtain the nice transformations we wanted, and now the issue is that we have to work on infinite-dimensional spaces.
\section{Group Action on the Infinite Prolongations}
\label{groupaction}

Now that we have our diffeomorphism between infinite jet bundles, we would like to know the form of our group action $\bfG$. Instead of working with a subgroup of $GL(n,\mathbb{R})$, what we have now is a group of transformations that take local coframings of $T^*\overline\Sigma_\infty$ to local coframings of $T^*\overline\Lambda_\infty$. In an equivalence problem of finite dimensional objects, $\varphi^* \bfo = \gamma \bfO$, $\gamma$ is essentially the pointwise Jacobian of the diffeomorphisms $\varphi$. The same is true in the case of infinite prolongations.

Suppose we have a transformation $(t,\bfx,\bfu,\dot{\bfu},\ldots)\mapsto(t,\bfy,\bfv,\dot{\bfv},\ldots)$ such that $t\mapsto t$. Suppose $\bfy=\bfy( \bfx,\bfu,\dot{\bfu},\ldots,\bfu^{(J)})$, i.e. $\bfy_{\bfu^{(J)}}$ is nonzero and $\bfy_{\bfu^{(k)}}=0$ for all $k>J$. If $\bfv=\bfv( \bfx,\bfu,\dot{\bfu},\ldots,\bfu^{(J_1)})$, we need to know first of all how $J_1$ is related to $J$.

On the one hand, we can directly compute the time derivative of $\bfy$ using the chain rule.
\begin{eqnarray*}
\frac{d\bfy}{dt} & = & \frac{d}{dt}\bfy( \bfx,\bfu,\dot{\bfu},\ldots,\bfu^{(J)} ) \\
 & = & \bfy_\bfx( \bfx,\bfu,\dot{\bfu},\ldots,\bfu^{(J)} ) \bff(\bfx,\bfu) + \bfy_\bfu( \bfx,\bfu,\dot{\bfu},\ldots,\bfu^{(J)} ) \dot{\bfu} \\
 & & \qquad +\ \ldots\ \bfy_{\bfu^{(J)}}( \bfx,\bfu,\dot{\bfu},\ldots,\bfu^{(J)} ) \bfu^{(J+1)}
\end{eqnarray*}
On the other hand, $\dot{\bfy}=\bfg(\bfy,\bfv)$.
\begin{eqnarray*}
\frac{d\bfy}{dt} & = & \bfg\left( \bfy ( \bfx,\bfu,\dot{\bfu},\ldots,\bfu^{(J)} ), \bfv ( \bfx,\bfu,\dot{\bfu},\ldots,\bfu^{(J_1)}) \right)
\end{eqnarray*}
Comparing these two versions of $\frac{d\bfy}{dt}$ shows that $\bfv=\bfv( \bfx,\bfu,\dot{\bfu},\ldots,\bfu^{(J+1)})$. Thus we have the following theorem.

\begin{theorem}
\label{nonzero}
$\bfy_{\bfu^{(J)}}$ is nonzero and $\bfy_{\bfu^{(k)}}=0$ for all $k>J$ if and only if $\bfv_{\bfu^{(J+1)}}$ is nonzero and $\bfv_{\bfu^{(k)}}=0$ for all $k>J+1$.
\end{theorem}
This relation and its repeated derivatives with respect to $t$ show that
\[
\bfv^{(i)}=\bfv^{(i)}( \bfx,\bfu,\dot{\bfu},\ldots,\bfu^{(J+i+1)}).
\]

Theorem \ref{nonzero} relates to our coframing as follows. Here we are omitting the pullbacks from our notation.
\begin{eqnarray*}
d\bfy & = & d\left( \bfy( \bfx,\bfu,\dot{\bfu},\ldots,\bfu^{(J)} ) \right) \\
 & = & \frac{\partial \bfy}{\partial \bfx}d\bfx + \sum_{i=0}^J \frac{\partial \bfy}{\partial \bfu^{(i)}} d\bfu^{(i)} \\
d\bfy - \bfg(\bfy,\bfv)dt & = & \frac{\partial \bfy}{\partial \bfx}d\bfx + \sum_{i=0}^J \frac{\partial \bfy}{\partial \bfu^{(i)}} d\bfu^{(i)} - \bfg\left( \bfx,\bfu,\ldots,\bfu^{(J+1)} \right)dt \\
 & = & \bA^0_0(d\bfx-\bff(\bfx,\bfu)dt) + \sum_{i=0}^J \bA^0_{i+1}\left( d\bfu^{(i)}-\bfu^{(i+1)}dt \right)
\end{eqnarray*}
where $\bA^0_j$, $0\leq j\leq J+1$, are matrices of functions of $\bfx,\bfu,\ldots,\bfu^{(J+1)}$. The fact that dynamic equivalence is time independent and takes solutions to solutions implies that there is no additional $\bA^0_{-1}\ dt$ here.

Similar calculations for $d\bfv^{(i)} - \bfv^{(i+1)}dt$ imply that our preferred coframings \eqref{coframe} transform in the following way,
\begin{eqnarray}
\label{Apullback}
\Phi_\infty^*\bfO = \bA\bfo & \qquad & \left({\Phi_\infty}^{-1}\right)^*\bfo = \left( {\Phi_\infty}^{-1} \bA \right)^{-1}\bfO
\end{eqnarray}
where $\bfO$,$\bA$,$\bfo$ have the form
\begin{eqnarray*}
\bfO = 
\left( \begin{array}{c}
	\Omega^{-1} \\
	\bfO^0 \\
	\bfO^1 \\
	\bfO^2 \\
	\vdots
\end{array}\right) , & &
\bfo = 
\left( \begin{array}{c}
	\omega^{-1} \\
	\bfo^0 \\
	\bfo^1 \\
	\bfo^2 \\
	\vdots \\
	\bfo^{J+1} \\
	\bfo^{J+2} \\
	\vdots
\end{array}\right) ,
\end{eqnarray*}
\[
\bA = 
\left( \begin{array}{cccccccccc}
1 & 0_{1\times n} & 0_{1\times s} & 0_{1\times s} & \cdots & 0_{1\times s} & 0_{1\times s} & 0_{1\times s} & 0_{1\times s} & \cdots  \\
0_{n\times 1} & \bA^0_0 & \bA^0_1 & \bA^0_2 & \cdots & \bA^0_{J+1} & 0_{s\times s} & 0_{s\times s} & 0_{s\times s} & \cdots \\
0_{s\times 1} & \bA_0^1 & \bA^1_1 & \bA^1_2 & \cdots & \bA^1_{J+1} & \bA^1_{J+2} & 0_{s\times s} & 0_{s\times s} & \cdots \\
0_{s\times 1} & \bA_0^2 & \bA^2_1 & \bA^2_2 & \cdots & \bA^2_{J+1} & \bA^2_{J+2} & \bA^2_{J+3} & 0_{s\times s} & \cdots \\
 & & & & \vdots \\
\end{array} \right) ,
\]
and the $\bA^i_j$ are submatrices of the following sizes.
\begin{center}
\begin{tabular}{|c|ccccc|}
\hline
matrix & $\bA^0_0$ & $\bA^0_j$ & $\bA^i_0$ & $\bA_j^i$ & ($i,j\geq 1$) \\
\hline
size & $n\times n$ & $n\times s$ & $s\times n$ & $s\times s$ & \\
\hline
\end{tabular}
\end{center}
A matrix $\bA$ of the above form for a fixed $J$ may have an inverse matrix similar to the above form with arbitrarily large $K$. For example, composition of dynamic equivalence maps leads to arbitrarily large $J$ and $K$.

From here on out, for any statement or theorem about $\bA$, an analogous statement or theorem also holds for $\bA^{-1}$ unless otherwise noted. These have been omitted for brevity. Submatrices of $\bA$ ($\bA^{-1}$) will be denoted by uppercase $\bA^i_j$ (lowercase $\ba^i_j$). Individual entries of these submatrices will denoted by $(A^i_j)^k_l$ ($(a^i_j)^k_l$), which are functions and thus not bolded. If a particular submatrix is in fact a scalar, which happens when $s=1$, then no bold face type will be used: $A^i_j$.

\begin{theorem}
\label{Arepeats}
Given a dynamic equivalence $\Phi_\infty^* \bfO = \bA \bfo$ with adapted coframings \eqref{coframe}, $\bA^i_{J+i+1} = \bA^1_{J+2}$ for all $i\geq 1$.
\end{theorem}

\begin{proof}
This proof is by induction on $i$. The case of $i=1$ is obvious. For $i\geq 1$, consider $d(\bfO^{i})$. Where an equivalence sign $\equiv$ is present below, it is because we are considering the equation modulo the linear span of $\{ \bfo^0,\ldots,\bfo^{J+i+1} \}$. Keep in mind that we are working with vector equations here. Recall that $\bfo^0 = d\bfx - \bff(\bfx,\bfu)dt$ is $n\times 1$, and $\bfo^{j} = d\bfu^{(j-1)} - \bfu^{(j)}dt$, $\bfO^{j} = d\bfv^{(j-1)} - \bfv^{(j)}dt$ are $s\times 1$ for $j\geq 1$. It is straightforward to verify in coordinates that $d\bfo^j=-\bfo^{j+1}\wedge dt$ and $d\bfO^j=-\bfO^{j+1}\wedge dt$ for $j\geq 1$.

On the one hand,
\begin{eqnarray*}
d(\bfO^{i}) & = & d\left(d\bfv^{(i-1)} - \bfv^{(i)}dt\right) \\
 & = & -d\bfv^{(i)}\wedge dt \\
 & = & -\bfO^{i+1}\wedge dt \\
 & \equiv & -\bA^{i+1}_{J+i+2} \bfo^{J+i+2}\wedge dt.
\end{eqnarray*}
On the other hand,
\begin{eqnarray*}
d(\bfO^{i}) & = & d( \sum_{j=0}^{J+i+1} \bA^{i}_{j} \bfo^{j} ) \\
 & = & \sum_{j=0}^{J+i+1} \left[ d(\bA^{i}_{j})\wedge\bfo^{j} + \bA^{i}_{j} d(\bfo^{j}) \right] \\
 & = & \sum_{j=0}^{J+i+1} \left[ d(\bA^{i}_{j})\wedge\bfo^{j} - \bA^{i}_{j} \bfo^{j+1}\wedge dt \right] \\
 & \equiv & - \bA^{i}_{J+i+1} \bfo^{J+i+2}\wedge dt.
\end{eqnarray*}
Since the $\bfo^j$ form a coframing, they are linearly independent. Thus we can conclude that
\begin{eqnarray*}
\bA^{i}_{J+i+1} & = & \bA^{i+1}_{J+i+2}.
\end{eqnarray*}
\end{proof}

While this does not completely characterize the group action of dynamic equivalence, it will be sufficient to prove a result in the next section that classifies dynamic equivalence in the case of one control variable. Later sections will narrow down what this group $\bA$ looks like; however, we will never completely characterize it. What we do prove about $\bA$ will be sufficient for some non-existence results.
\section{Scalar Control}
\label{scalarcontrol}

The following theorem about dynamic equivalence in the case of one control variable has been known for some time. What is presented here is a proof based on Pomet's work \cite{po1} that has been adapted to this framework of coframings on infinite jet bundles. It reduces all dynamic equivalences of control systems with just one control variable to the case of static equivalence.

\begin{theorem}
Let the control systems $\Sigma$, $\Lambda$ in \eqref{cs_eq} be dynamically equivalent with $s=1$ control variable and $m$, $n$ state variables, respectively. If $m=n$, then the systems are in fact static equivalent. If $m<n$, then the systems are static equivalent after a finite number of prolongations of the smaller system $\Sigma$.
\end{theorem}

\begin{proof}
Let $\bA=(\bA^i_j)$ and $\bA^{-1}=(\ba^i_j)$ as before.

If $m<n$, prolong $\Sigma$ until $m=n$. Suppose the coframings of $\Sigma$, $\Lambda$ in \eqref{cs_eq} pull back as in \eqref{Apullback}. Suppose there exist nonnegative integers $J$ and $K$ such that $\bfx_{v^{(J)}}$ and $\bfy_{u^{(K)}}$ are nonzero. In Theorem \ref{static} in the next section, it is shown that it is not possible for just one of $J$ or $K$ to be $-1$, i.e. $\bA^i_j={\bf 0}$ for all $j>i$ if and only if $\ba^i_j={\bf 0}$ for all $j>i$. So both $J$ and $K$ must be nonnegative for a strict dynamic equivalence to exist.

By the computations in the previous section, $\bA^0_{J+1}$ is a nonzero $n\times 1$ matrix. Likewise, $a^i_{K+i+1}$ is a nonzero function for all $i\geq 1$. Because $\bA^0_{J+1}$ is a nonzero $n\times 1$ vector, and $a^{J+1}_{K+J+2}$ is a nonzero function, their product $\bA^0_{J+1} a^{J+1}_{K+J+2}$ is a nonzero $n\times 1$ vector. However $\bA \bA^{-1}$ is the identity. Therefore $\bA^0_{J+1} a^{J+1}_{K+J+2}$, which is an off diagonal $n\times 1$ entry since $0<K+J+2$, must be an all zero $n\times 1$ vector.

This is a contradiction. Thus $J$ and $K$ cannot exist, and $\bfx_{v^{(J)}} = \bfy_{u^{(K)}} = 0$ for all $J,K\geq 0$. This shows that the equivalence is in fact static.
\end{proof}

\section{Group Adaptations for Two Controls}
\label{2controls}

The last section dealt with the case of a scalar control, in which dynamic and static equivalence are one and the same. Now we will work on the next simplest case of two controls ($s=2$) with $J=K=0$. In the case of one control variable, there is essentially no ``room for freedom" to allow a true dynamic equivalence, aside from prolongations. With two control variables, there is now ``room" to have a strict dynamic equivalence, but just barely. While larger values of $J$ and $K$ increase the flexibility of possible dynamic equivalences, in this section we will show that there is really only one way to have a strict dynamic equivalence of two systems with $J=K=0$.

\subsection{Nonautonomous Static Equivalence}

Recall the notation we have developed thus far for the pullbacks of our preferred coframings. Note the equivalent submatrices $\bA^i_{i+1}$, $i\geq 1$, from Theorem \ref{Arepeats}.

\[
\bar\Phi_\infty^*
\left( \begin{array}{c}
\Omega^{-1} \\
\bfO^0 \\
\bfO^{1} \\
\bfO^{2} \\
\vdots
\end{array} \right) = 
\left( \begin{array}{ccccccc}
1 & 0_{1\times n} & 0_{1\times 2} & 0_{1\times 2} & 0_{1\times 2} & 0_{1\times 2} & \cdots  \\
0_{n\times 1} & { \bA^0_0} & { \bA^0_1} & 0_{n\times 2} & 0_{n\times 2} & 0_{n\times 2} & \cdots  \\
0_{2\times 1} & { \bA_0^1} & { \bA^1_1} & { \bA^1_2} & 0_{2\times 2} & 0_{2\times 2} & \cdots \\
0_{2\times 1} & { \bA_0^2} & { \bA^2_1} & { \bA^2_2} & { \bA^1_{2}} & 0_{2\times 2} & \cdots \\
 & & & & \vdots \\
\end{array} \right)
\left( \begin{array}{c}
\omega^{-1} \\
\bfo^0 \\
\bfo^{1} \\
\bfo^{2} \\
\bfo^{3} \\
\vdots
\end{array} \right)
\]

\[
\left(\bar\Phi_\infty^{-1}\right)^*
\left( \begin{array}{c}
\omega^{-1} \\
\bfo^0 \\
\bfo^{1} \\
\bfo^{2} \\
\vdots
\end{array} \right) = 
\left( \begin{array}{ccccccc}
1 & 0_{1\times n} & 0_{1\times 2} & 0_{1\times 2} & 0_{1\times 2} & 0_{1\times 2} & \cdots  \\
0_{n\times 1} & { \ba^0_0} & { \ba^0_1} & 0_{n\times 2} & 0_{n\times 2} & 0_{n\times 2} & \cdots  \\
0_{2\times 1} & { \ba_0^1} & { \ba^1_1} & { \ba^1_2} & 0_{2\times 2} & 0_{2\times 2} & \cdots \\
0_{2\times 1} & { \ba_0^2} & { \ba^2_1} & { \ba^2_2} & { \ba^1_{2}} & 0_{2\times 2} & \cdots \\
 & & & & \vdots \\
\end{array} \right)
\left( \begin{array}{c}
\Omega^{-1} \\
\bfO^0 \\
\bfO^{1} \\
\bfO^{2} \\
\bfO^{3} \\
\vdots
\end{array} \right)
\]

In what follows, we will refer to a group element $\bfgg$,
\[
\bfgg = \left( \begin{array}{ccccccc}
1 & 0_{1\times n} & 0_{1\times 2} & 0_{1\times 2} & 0_{1\times 2} & 0_{1\times 2} & \cdots  \\
0_{n\times 1} & { \bfg^0_0} & 0_{n\times 2} & 0_{n\times 2} & 0_{n\times 2} & 0_{n\times 2} & \cdots  \\
0_{2\times 1} & { \bfg_0^1} & { \bfg^1_1} & 0_{2\times 2} & 0_{2\times 2} & 0_{2\times 2} & \cdots \\
0_{2\times 1} & { \bfg_0^2} & { \bfg^2_1} & { \bfg^2_2} & 0_{2\times 2} & 0_{2\times 2} & \cdots \\
 & & & & \vdots \\
\end{array} \right),
\]
that acts on our coframings as \textbf{nonautonomous static equivalence}, meaning $\bfgg^i_j={\bf 0}$ for all $i<j$. This terminology arises from the fact that such $\bfgg$ arise as the Jacobian of a time-dependent static equivalence $\tilde\bfx=\bfx(\bfx,t)$ on the contact system of the infinite prolongation $\overline{\Sigma}_\infty$. Unlike the matrix representing a true static equivalence, $\bfgg$ allows changes of variables such as $x_i\mapsto x_i+t$. Note that such equivalences take a coframing on $\overline\Sigma_\infty$ to another coframing on $\overline\Sigma_\infty$ (or $\overline\Lambda_\infty$ to $\overline\Lambda_\infty$).

As in the case of dynamic equivalence, we wish to require that the following structure equations are preserved by nonautonomous static equivalence.
\begin{eqnarray}
\label{nice_struct_eqns}
d\bfO^i \in \textrm{span }\{ \bfO^{i+1}\wedge\bfO^{-1} \} & \textrm{ mod } \bfO^j, & \ 0\leq j \leq i.
\end{eqnarray}
This additional condition allows us to simplify the form of $\bfgg$ much like we did for $\bA$ in Theorem \ref{Arepeats}. The proof is identical to that of Theorem \ref{Arepeats} with $J=-1$.

\begin{theorem}
Given a nonautonomous static equivalence $\bfgg$ that preserves the structure equations \eqref{nice_struct_eqns}, $\bfg^i_i = \bfg^1_1$ for all $i\geq 1$.
\end{theorem}

A straightforward calculation in coordinates shows that every static equivalence is a nonautonomous static equivalence, but of course the converse is not true.

Later we will be showing that every dynamic equivalence with $J=K=0$ can be factored into a constant matrix composed with nonautonomous static equivalences. This result will be key in proving the main classification results of this paper.

\begin{theorem}
\label{static}
$\bA$ is a nonautonomous static equivalence, i.e. $\bA^i_j={\bf 0}$ for all $i<j$, if and only if $\bA^{-1}$ is also a nonautonomous static equivalence.
\end{theorem}

\begin{proof}
If $\bA^0_1={\bf 0}$, then $\bA^0_0$ is a rank $n$ matrix, hence invertible. Let the submatrices of $\bA^{-1}$ be denoted by $\ba^i_j$. Since $\bA\bA^{-1}=Id$, the off diagonal element $\bA^0_0\ba^0_1$ must be zero. Because $\bA^0_0$ is invertible, this means $\ba^0_1=0$. By Theorem \ref{nonzero}, $\ba^0_1$ is zero if and only if $\ba^1_2$ is too. Theorem \ref{Arepeats} completes the proof since $\ba^i_{i+1} = \ba^1_2$ for all $i\geq 1$ and $\ba^i_j=0$ for all $j>i$. Therefore $\bA^{-1}$ is a nonautonomous static equivalence.
\end{proof}

\subsection{Factoring $\bA$}

In the following section, we will prove several theorems about the rank of certain submatrices of $\bA$. This chapter will culminate in the final theorem, theorem \eqref{factorA}, which states that we can factor our dynamic equivalence in a special way: $\bA=\bfgg\bfS\bfG$. The $\bfgg$ and $\bfG$ are two nonautonomous static equivalences which encapsulate the traditional change of variables, as in static equivalence. The $\bfS$ is a fixed constant orthogonal matrix which incorporates the mixing of higher derivatives into dynamic equivalence.

\begin{theorem}
Given a strictly dynamic equivalence $\bA$ with $s=2$ and $J=K=0$, $\bA^1_2$ (a $2\times 2$ submatrix) has rank 1.
\end{theorem}

\begin{proof}
We know that $\bA^1_2$ cannot have rank zero by Theorem \ref{static}. Assume the rank of $\bA^1_2$ is two. Then through a change of coframing $\tilde\bfo = \bfG\bfo$ via static equivalence $\bfG$, it can be arranged that the elements of $\tilde \bA = \bA \bfG^{-1}$ look as follows.
\[
\tilde\Phi_\infty^*
\left( \begin{array}{@{}c@{}}
\Omega^{-1} \\
\bfO^0 \\
\bfO^{1} \\
\bfO^{2} \\
\vdots
\end{array} \right) =
\left( \begin{array}{ccccccc}
1 & 0_{1\times n} & 0_{1\times 2} & 0_{1\times 2} & 0_{1\times 2} & 0_{1\times 2} & \cdots  \\
0_{n\times 1} & { \tilde \bA^0_0} & { \tilde \bA^0_1} & 0_{n\times 2} & 0_{n\times 2} & 0_{n\times 2} & \cdots  \\
0_{2\times 1} & 0_{2\times 2} & 0_{2\times 2} & Id_{2\times 2} & 0_{2\times 2} & 0_{2\times 2} & \cdots \\
0_{2\times 1} & 0_{2\times 2} & 0_{2\times 2} & 0_{2\times 2} & Id_{2\times 2} & 0_{2\times 2} & \cdots \\
 & & & & \vdots \\
\end{array} \right)
\left( \begin{array}{@{}c@{}}
\tilde\omega^{-1} \\
\tilde\bfo^0 \\
\tilde\bfo^{1} \\
\tilde\bfo^{2} \\
\tilde\bfo^{3} \\
\vdots
\end{array} \right)
\]
We have $\tilde\Phi_\infty^*\bfO^j = \tilde\bfo^{j+1}$ for $j\leq 1$. By the nature of pullbacks, this also means $\left(\tilde\Phi_\infty^{-1}\right)^*\tilde\bfo^{j+1} = \bfO^j$. However this means that $\bA^{-1}$ now looks as follows.
\[
\left(\tilde\Phi_\infty^{-1}\right)^*
\left( \begin{array}{@{}c@{}}
\tilde\omega^{-1} \\
\tilde\bfo^0 \\
\tilde\bfo^{1} \\
\tilde\bfo^{2} \\
\vdots
\end{array} \right) = 
\left( \begin{array}{ccccccc}
1 & 0_{1\times n} & 0_{1\times 2} & 0_{1\times 2} & 0_{1\times 2} & 0_{1\times 2} & \cdots  \\
0_{n\times 1} & { \tilde \ba^0_0} & { \tilde \ba^0_1} & 0_{n\times 2} & 0_{n\times 2} & 0_{n\times 2} & \cdots  \\
0_{2\times 1} & { \tilde \ba_0^1} & { \tilde \ba^1_1} & { \tilde \ba^1_2} & 0_{2\times 2} & 0_{2\times 2} & \cdots \\
0_{2\times 1} & 0_{2\times 2} & Id_{2\times 2} & 0_{2\times 2} & 0_{2\times 2} & 0_{2\times 2} & \cdots \\
 & & & & \vdots \\
\end{array} \right)
\left( \begin{array}{@{}c@{}}
\Omega^{-1} \\
\bfO^0 \\
\bfO^{1} \\
\bfO^{2} \\
\bfO^{3} \\
\vdots
\end{array} \right)
\]
In particular, $0=\tilde\ba^2_3=\tilde\ba^1_2$. By the above argument this means that $\tilde\ba^0_1=0$ and the equivalence is static. This contradicts $J=K=0$. Therefore the rank of $\bA^1_2$ must be one.
\end{proof}

\begin{theorem}
Given a dynamic equivalence $\bA$ with $s=2$ and $J=K=0$, $\bA^0_1$ (an $n\times 2$ submatrix) has rank 1.
\end{theorem}

\begin{proof}
The rank of $\bA^0_1$ is either 0, 1, or 2. If the rank is zero, then the equivalence is static. Consider $(\bA\bA^{-1})^0_2 = \bA^0_1\ba^1_2 = 0_{n\times 2}$. If rank of $\bA^0_1$ is two, then $\bA^0_1$ has a $2\times n$ left inverse, and we conclude $\ba^1_2=0_{2\times 2}$. However this again implies a static equivalence, so the rank is not two.
\end{proof}

The plan now is to use this knowledge of the ranks to normalize $\bA$ via non-autonomous static group actions to $\bfo$ and $\bfO$. This will isolate the dynamic part of the mapping to one very specific form $\bfbO = \bfS \bfbo$, where $\bfbO = \bfgg^{-1}\bfO$, $\bfbo=\bfG\bfo$, $\bfo$ and $\bfO$ are our preferred coframings (to be determined later), and $\bfgg$, $\bfG$ are non-autonomous static group elements. An explicit example of how this is done will follow in the next section.

Starting with the fact that $\bA^0_1$ has rank one, we know it can be normalized to the following form through Gauss-Jordan elimination, which in this context is nonautonomous static equivalences applied to the coframings $\bfo$ and $\bfO$.
\[
\bA^0_1 = \left(\begin{array}{cc}
0 & 1 \\
0 & 0 \\
\vdots & \vdots \\
0 & 0 \\
\end{array}\right)
\]
Recall that all $\bfo^i = (\omega^i_j)$ and $\bfO^i = (\Omega^i_j)$ are vectors. If we add multiples of $\omega^0_i$ to $\omega^1_1$, we can eliminate the first row of $\bA^0_0$. Note that this can be accomplished by a static group action.
\[
\bA^0_0 = \left(\begin{array}{ccc}
0 & \cdots & 0 \\
* & \cdots & * \\
\vdots & & \vdots \\
* & \cdots & * \\
\end{array}\right)
\]
Since the $n\times (n+2)$ matrix $(\bA^0_0\ | \bA^0_1)$ must have rank $n$ for $\bA$ to be invertible, the last $n-1$ rows of $\bA^0_0$ must have rank $n-1$. This allows us to normalize the rest of $\bA^0_0$ via a static group action.
\[
\bA^0_0 = \left(\begin{array}{cc}
0 & 0_{1\times (n-1)} \\
0 & Id_{(n-1)\times (n-1)} \\
\end{array}\right)
\]
The first $n+1$ rows of $\bA$ have now been reduced to ones and zeros.

Since the rank of $\bA^i_{i+1}$ is one, there exists a  non-autonomous static equivalences applied to both coframings $\bfo$ and $\bfO$ that yields a new coframing with
\[
\bA^i_{i+1} = \left(\begin{array}{cc}
0 & 0 \\
0 & 1
\end{array}\right),
\]
Everything to the left of the ones in each $\bA^i_{i+1}$ can be eliminated by a nonautonomous static equivalence that redefines $\omega^i_{i+1}$. In fact anything to the left of or below a one in the matrix $\bA$ can essentially be absorbed by a non-autonomous static equivalence that redefines either $\bfbo$ (horizontal zeros) or $\bfbO$ (vertical zeros). For ease of notation, these newly redefined coframings, which differ from the original preferred coframings by non-autonomous static equivalences, will be still be denoted with $\bfO$ and $\bfo$. This leaves the following simplified form of $\bA$.
{\small
\[
\left( \begin{array}{c|cc|cc|cc|cc|ccc}
1 & 0 & {\bf 0} & 0 & 0 & 0 & 0 & 0 & 0 & 0 & 0 & \cdots  \\
\hline
0 & 0 & {\bf 0} & 0 & 1 & 0 & 0 & 0 & 0 & 0 & 0 & \cdots \\
0 & {\bf 0} & Id_{(n-1)\times(n-1)} & {\bf 0} & {\bf 0} & {\bf 0} & {\bf 0} & {\bf 0} & {\bf 0} & {\bf 0} & {\bf 0} & \cdots  \\
\hline
0 & (\bA^1_0)^1_1 & {\bf 0} & (\bA^1_1)^1_1 & 0 & 0 & 0 & 0 & 0 & 0 & 0 & \cdots \\
0 & 0 & {\bf 0} & 0 & 0 & 0 & 1 & 0 & 0 & 0 & 0 & \cdots \\
\hline
0 & (\bA^2_0)^1_1 & {\bf 0} & (\bA^2_1)^1_1 & 0 & (\bA^2_2)^1_1 & 0 & 0 & 0 & 0 & 0 & \cdots \\
0 & 0 & {\bf 0} & 0 & 0 & 0 & 0 & 0 & 1 & 0 & 0 & \cdots \\
\hline
0 & (\bA^3_0)^1_1 & {\bf 0} & (\bA^3_1)^1_1 & 0 & (\bA^3_2)^1_1 & 0 & (\bA^3_3)^1_1 & 0 & 0 & 0 & \cdots \\
0 & 0 & {\bf 0} & 0 & 0 & 0 & 0 & 0 & 0 & 0 & 0 & \cdots \\
 & & & & \vdots & & & & & \\
\end{array} \right)
\]
}

Now if $(\bA^1_0)^1_1$ is zero, one of the other $(\bA^i_0)^1_1$, $i>1$, must be nonzero. This follows from the fact that $\bA^{-1}\bA = {\bf Id}$, in particular $((\bA^{-1}\bA)^0_0)^1_1 = 1$. If $(\bA^1_0)^1_1$ is zero, there is an $i>1$ such that $(\ba^0_i)^1_1 (\bA^i_0)^1_1$ is nonzero. But $(\ba^0_i)^1_1$ being nonzero implies $K+1\geq i>1$. Since we are restricting our consideration to $K=0$, this cannot happen. Therefore $(\bA^1_0)^1_1$ must be nonzero. Since $(\bA^1_0)^1_1$ is nonzero, it can be scaled to unity through a nonautonomous static group action. All of the other $(\bA^i_0)^1_1$ can then be eliminated through non-autonomous static group actions (adding multiples of rows in this case).

It can similarly be shown that when $J=K=0$, $(\bA^{i+1}_i)^1_1$ is nonzero and can be scaled to unity. All entries below them can be made zero. By examining $\bA^{-1}\bA = Id$ one can also check that any of the $(\bA^i_i)^1_1$ being nonzero leads to $K+1\geq 2$, and therefore $(\bA^i_i)^1_1=0$.

Finally all the group freedom of $\bA$ has been utilized through non-autonomous static group actions on $\bfo$ and $\bfO$, and what is left is the following constant matrix.

\begin{multline}
\label{stackpole}
\bfS = \\
{\textrm{\footnotesize
$\left( \begin{array}{c|c|c@{\hspace{1pt}}c|cc|cc|c}
1 & 0_{1\times n} & 0 & 0 & 0 & 0 & 0 & 0 & \cdots  \\
\hline
0_{n\times 1} \! &\!
	\begin{array}{cc}
	0 & 0_{1\times (n-1)} \\
	0 \! &\! Id_{(n-1)\times (n-1)} \\
	\end{array}
 & \begin{array}{c}
	0 \\
	0_{(n-1)\times 1} \\
	\end{array} \! &\! 
	\begin{array}{c}
	1 \\
	0_{(n-1)\times 1} \\
	\end{array} \! &\! 0_{n\times 1} \! &\! 0_{n\times 1} \! &\! 0_{n\times 1} \! &\! 0_{n\times 1} \! &  \cdots  \\
\hline
0 & \begin{array}{cc}
1 & 
0_{1\times (n-1)} \\
\end{array} & 0 & 0 & 0 & 0 & 0 & 0 & \cdots \\
0 & \begin{array}{cc}
0 & 
0_{1\times (n-1)} \\
\end{array} & 0 & 0 & 0 & 1 & 0 & 0 & \cdots \\
\hline
0 & 0_{1\times n} & 1 & 0 & 0 & 0 & 0 & 0 & \cdots \\
0 & 0_{1\times n} & 0 & 0 & 0 & 0 & 0 & 1 & \cdots \\
\hline
 & & & & \vdots & & & \\
\end{array} \right)
$}}
\end{multline}

It is easy to check that $\bfS$ is orthogonal, i.e. $\bfS^{-1}=\bfS^T$. Thus we have proved the following theorem.

\begin{theorem}
\label{factorA}
Given preferred coframings $\bfo$ and $\bfO$ \eqref{coframe} on $\overline{\Sigma}_\infty$ and $\overline{\Lambda}_\infty$ for control systems $\Sigma$ and $\Lambda$, respectively, with $s=2$ and a dynamic equivalence $\Phi_\infty$ with $J=K=0$ taking $\Sigma_\infty$ to $\Lambda_\infty$, the coframing pulls back as follows:
\begin{eqnarray*}
\bar\Phi_\infty^* \bfO & = & \bfgg\ \bfS\ \bfG\ \bfo
\end{eqnarray*}
where $\bfgg$ and $\bfG$ are nonautonomous static equivalences and $\bfS$ is given by \eqref{stackpole} above.
\end{theorem}

This theorem means that, up to nonautonomous static equivalence, a dynamic equivalence with $J=K=0$ has a very specific form which is encoded in this specific orthogonal matrix $\bfS$. Most of the apparent complexity of dynamic equivalence actually arises from static equivalence on either side, and the essence of dynamic equivalence is actually quite simple.
\section{Factoring the Dynamic Equivalence: An Example}
\label{factoring}

What we have shown so far is that, given a dynamic equivalence $\bfO = \bA\bfo$ where $J=K=0$ and $s=2$ ($n$ is still arbitrary), we can decompose the group action $\bA = \bfgg \bfS \bfG$ where $\bfS$ is defined by \eqref{stackpole} and $\bfG$ and $\bfgg$ are non-autonomous static equivalent group elements, i.e.
\begin{eqnarray*}
\bfG & = & 
\left( \begin{array}{ccccccc}
1 & 0_{1\times n} & 0_{1\times 2} & 0_{1\times 2} & 0_{1\times 2} & 0_{1\times 2} & \cdots  \\
0_{n\times 1} & { \bfG^0_0} & 0_{n\times 2} & 0_{n\times 2} & 0_{n\times 2} & 0_{n\times 2} & \cdots  \\
0_{2\times 1} & { \bfG_0^1} & { \bfG^1_1} & 0_{2\times 2} & 0_{2\times 2} & 0_{2\times 2} & \cdots \\
0_{2\times 1} & { \bfG_0^2} & { \bfG^2_1} & { \bfG^2_2} & 0_{2\times 2} & 0_{2\times 2} & \cdots \\
 & & & & \vdots \\
\end{array} \right)
\end{eqnarray*}
Recall that nonautonomous static equivalence is not a true static equivalence
\[
(\bfx,\bfu) \mapsto (\bfy(\bfx), \bfv(\bfx,\bfu)).
\]
Unlike static equivalence, which is autonomous (time-independent), nonautonomous static equivalence can have explicit time dependence, for example, $x_i \mapsto x_i + t$. Its group action does not preserve the algebraic ideal $\{d\bfx\}$, just the ideal $\{d\bfx - \bff(\bfx,\bfu)\ dt\}$. This equivalence is more general than static equivalence.

Let us phrase the problem now as follows. Dynamic equivalence looks like $\bfO = \bA\bfo$ where $\bA = \bfgg \bfS \bfG$. We can attack this problem in steps. First we will consider the coframing $\bfbO = \bfS \bfG \bfo$. Then what remains will be the non-autonomous static problem $\bfO = g\bfbO$.

\begin{example}\end{example}

Let us consider the following two dynamically equivalent systems:
\begin{eqnarray*}
\dot{x}_1 = u_1 & \qquad & \dot{y}_1 = v_1 \\
\dot{x}_2 = u_2 & \qquad & \dot{y}_2 = v_2 \\
\dot{x}_3 = x_2u_1 & \qquad & \dot{y}_3 = y_2
\end{eqnarray*}
An actual dynamic equivalence is given by the following maps between the infinite jet bundles.
\begin{eqnarray*}
\Phi_\infty( \bfx,\bfu,\dot{\bfu},\ldots ) & = & (\ x_1x_2 - x_3,\ u_2,\ x_2,\ x_1u_2,\ \dot{u}_2,\ \ldots\ ) \\
\displaystyle\Phi_\infty^{-1}( \bfy,\bfv,\dot{\bfv},\ldots ) & = & (\ v_1/y_2,\ y_3,\ y_3v_1/y_2 - y_1,\ \ldots\ )
\end{eqnarray*}
Here is a coframing for each of the infinite jet bundles. The choice of $\omega^0_3$, while not obvious, is not arbitrary. We will see why in a later section. For this example only the $dt$ piece of the coframing has been left out. Since $t\mapsto t$, this would just add a one and many zeros to the matrices.
\[
\begin{array}{lcl}
\bfo^0 = \left(\begin{array}{c}
dx_1 - u_1 dt \\
dx_2 - u_2 dt \\
dx_3 - x_2u_1 dt - x_2(dx_1 - u_1 dt) \\
\end{array}\right) & \qquad\qquad & 
\bfO^0 = \left(\begin{array}{c}
dy_1 - v_1 dt \\
dy_2 - v_2 dt \\
dy_3 - y_2 dt
\end{array}\right) \\
\bfo^1 = \left(\begin{array}{c}
du_1 - \dot{u}_1 dt \\
du_2 - \dot{u}_2 dt
\end{array}\right) & \qquad\qquad & 
\bfO^1 = \left(\begin{array}{c}
dv_1 - \dot{v}_1 dt \\
dv_2 - \dot{v}_2 dt
\end{array}\right) \\
\bfo^2 = \left(\begin{array}{c}
d\dot{u}_1 - \ddot{u}_1 dt \\
d\dot{u}_2 - \ddot{u}_2 dt
\end{array}\right) & \qquad\qquad & 
\bfO^2 = \left(\begin{array}{c}
d\dot{v}_1 - \ddot{v}_1 dt \\
d\dot{v}_2 - \ddot{v}_2 dt
\end{array}\right) \\
\vdots & & \vdots
\end{array}
\]

The pullback of $\bar\Phi_\infty$ is straightforward to calculate. For the rest of this section the pullback notation will be suppressed in order to emphasize and clarify the methods being used.
\begin{eqnarray*}
\Omega^0_1 & = & dy_1 - v_1 dt \\
& = & d(x_1x_2 - x_3) - (x_1u_2)\ dt \\
& = & x_2\ dx_1 + x_1\ dx_2 - dx_3 - (x_1u_2)\ dt - (x_2u_1)\ dt + (x_2u_1)\ dt \\
& = & -\left[ dx_3 - x_2u_1 dt - x_2(dx_1 - u_1 dt) \right] + x_1 \left( dx_2 - u_2\ dt \right) \\
& = & -\omega^0_3 + x_1\omega^0_2 \\
\\
\Omega^0_2 & = & dy_2 - v_2 dt \\
& = & du_2 - \dot{u}_2\ dt \\
& = & \omega^1_2 \\
\\
\Omega^0_3 & = & dy_3 - y_2 dt \\
& = & dx_2 - u_2 dt \\
& = & \omega^0_2 \\
\end{eqnarray*}
\begin{eqnarray*}
\Omega^1_1 & = & dv_1 - \dot{v}_1\ dt \\
& = & d(x_1u_2) - (u_1u_2 + x_1\dot{u}_2)\ dt \\
& = & u_2\ \left( dx_1 - u_1\ dt \right) + x_1 \left( du_2 - \dot{u}_2\ dt \right) \\
& = & u_2\omega^0_1 + x_1\omega^1_2 \\
\\
\Omega^1_2 & = & dv_2 - \dot{v}_2\ dt \\
& = & d\dot{u}_2 - \ddot{u}_2\ dt \\
& = & \omega^2_2 \\
& \vdots
\end{eqnarray*}

The pullback put in matrix form looks as follows.
\begin{eqnarray*}
\left(\begin{array}{c}
\Omega^0_1 \\
\Omega^0_2 \\
\Omega^0_3 \\
\hline
\Omega^1_1 \\
\Omega^1_2 \\
\hline
\Omega^2_1 \\
\Omega^2_2 \\
\vdots
\end{array}\right) & = & \left(\begin{array}{ccc|cc|cc|ccc}
0 & x_1 & -1 & 0 & 0 & 0 & 0 & 0 & 0 & \cdots \\
0 & 0 & 0 & 0 & 1 & 0 & 0 & 0 & 0 & \cdots \\
0 & 1 & 0 & 0 & 0 & 0 & 0 & 0 & 0 & \cdots \\
\hline
u_2 & 0 & 0 & 0 & x_1 & 0 & 0 & 0 & 0 & \cdots \\
0 & 0 & 0 & 0 & 0 & 0 & 1 & 0 & 0 & \cdots \\
\hline
\dot{u}_2 & 0 & 0 & u_2 & u_1 & 0 & x_1 & 0 & 0 & \cdots \\
0 & 0 & 0 & 0 & 0 & 0 & 0 & 0 & 1 & \cdots \\
\vdots & \vdots & \vdots & \vdots & \vdots & \vdots & \vdots & \vdots & \vdots \\
\end{array}\right)
\left(\begin{array}{c}
\omega^0_1 \\
\omega^0_2 \\
\omega^0_3 \\
\hline
\omega^1_1 \\
\omega^1_2 \\
\hline
\omega^2_1 \\
\omega^2_2 \\
\hline
\omega^3_1 \\
\omega^3_2 \\
\vdots
\end{array}\right)
\end{eqnarray*}

We will now follow the algorithm for producing $\bfS$. This amounts to a series of row or column operations which are static equivalences on the $\bfO$ or $\bfo$ coframes respectively. We will use the notation of a typical introduction to linear algebra course to represent these operations, i.e. $R_2\rightarrow R_2+R_3$ means to replace row 2 with row 2 plus row 3. Note that not every row operation is a legal static equivalence. For example, $R_1\rightarrow R_1 + R_4$ amounts to $x\mapsto x+u$, which is dynamic, not static.

First perform the following operations:
\[
R_1 \rightarrow R_3 \rightarrow R_2 \rightarrow R_1
\]
which results in the following coframing.
\begin{eqnarray*}
\left(\begin{array}{c}
\Omega^0_2 \\
\Omega^0_3 \\
\Omega^0_1 \\
\hline
\Omega^1_1 \\
\Omega^1_2 \\
\hline
\Omega^2_1 \\
\Omega^2_2 \\
\vdots
\end{array}\right) & = & \left(\begin{array}{ccc|cc|cc|ccc}
0 & 0 & 0 & 0 & 1 & 0 & 0 & 0 & 0 & \cdots \\
0 & 1 & 0 & 0 & 0 & 0 & 0 & 0 & 0 & \cdots \\
0 & x_1 & -1 & 0 & 0 & 0 & 0 & 0 & 0 & \cdots \\
\hline
u_2 & 0 & 0 & 0 & x_1 & 0 & 0 & 0 & 0 & \cdots \\
0 & 0 & 0 & 0 & 0 & 0 & 1 & 0 & 0 & \cdots \\
\hline
\dot{u}_2 & 0 & 0 & u_2 & u_1 & 0 & x_1 & 0 & 0 & \cdots \\
0 & 0 & 0 & 0 & 0 & 0 & 0 & 0 & 1 & \cdots \\
\vdots & \vdots & \vdots & \vdots & \vdots & \vdots & \vdots & \vdots & \vdots \\
\end{array}\right)
\left(\begin{array}{c}
\omega^0_1 \\
\omega^0_2 \\
\omega^0_3 \\
\hline
\omega^1_1 \\
\omega^1_2 \\
\hline
\omega^2_1 \\
\omega^2_2 \\
\hline
\omega^3_1 \\
\omega^3_2 \\
\vdots
\end{array}\right)
\end{eqnarray*}
Next perform the operation
\[
R_3 \rightarrow x_1R_2 - R_3
\]
to get this new coframing.
\begin{eqnarray*}
\left(\begin{array}{c}
\Omega^0_2 \\
\Omega^0_3 \\
x_1\Omega^0_3 - \Omega^0_1 \\
\hline
\Omega^1_1 \\
\Omega^1_2 \\
\hline
\Omega^2_1 \\
\Omega^2_2 \\
\vdots
\end{array}\right) & = & \left(\begin{array}{ccc|cc|cc|ccc}
0 & 0 & 0 & 0 & 1 & 0 & 0 & 0 & 0 & \cdots \\
0 & 1 & 0 & 0 & 0 & 0 & 0 & 0 & 0 & \cdots \\
0 & 0 & 1 & 0 & 0 & 0 & 0 & 0 & 0 & \cdots \\
\hline
u_2 & 0 & 0 & 0 & x_1 & 0 & 0 & 0 & 0 & \cdots \\
0 & 0 & 0 & 0 & 0 & 0 & 1 & 0 & 0 & \cdots \\
\hline
\dot{u}_2 & 0 & 0 & u_2 & u_1 & 0 & x_1 & 0 & 0 & \cdots \\
0 & 0 & 0 & 0 & 0 & 0 & 0 & 0 & 1 & \cdots \\
\vdots & \vdots & \vdots & \vdots & \vdots & \vdots & \vdots & \vdots & \vdots \\
\end{array}\right)
\left(\begin{array}{c}
\omega^0_1 \\
\omega^0_2 \\
\omega^0_3 \\
\hline
\omega^1_1 \\
\omega^1_2 \\
\hline
\omega^2_1 \\
\omega^2_2 \\
\hline
\omega^3_1 \\
\omega^3_2 \\
\vdots
\end{array}\right)
\end{eqnarray*}
The first three rows of the transformation now look like the first three rows of $\bfS$. Continue by letting
\[
R_4 \rightarrow \D \frac{R_4 - x_1R_1}{u_2}
\]
to yield the coframing below.
\begin{eqnarray*}
\left(\begin{array}{c}
\Omega^0_2 \\
\Omega^0_3 \\
x_1\Omega^0_3 - \Omega^0_1 \\
\hline
\left( \Omega^1_1 - x_1\Omega^0_2 \right) / u_2 \\
\Omega^1_2 \\
\hline
\Omega^2_1 \\
\Omega^2_2 \\
\vdots
\end{array}\right) & = & \left(\begin{array}{ccc|cc|cc|cc@{\hspace{3pt}}c}
0 & 0 & 0 & 0 & 1 & 0 & 0 & 0 & 0 & \cdots \\
0 & 1 & 0 & 0 & 0 & 0 & 0 & 0 & 0 & \cdots \\
0 & 0 & 1 & 0 & 0 & 0 & 0 & 0 & 0 & \cdots \\
\hline
1 & 0 & 0 & 0 & 0 & 0 & 0 & 0 & 0 & \cdots \\
0 & 0 & 0 & 0 & 0 & 0 & 1 & 0 & 0 & \cdots \\
\hline
\dot{u}_2 & 0 & 0 & u_2 & u_1 & 0 & x_1 & 0 & 0 & \cdots \\
0 & 0 & 0 & 0 & 0 & 0 & 0 & 0 & 1 & \cdots \\
\vdots & \vdots & \vdots & \vdots & \vdots & \vdots & \vdots & \vdots & \vdots \\
\end{array}\right)
\left(\begin{array}{c}
\omega^0_1 \\
\omega^0_2 \\
\omega^0_3 \\
\hline
\omega^1_1 \\
\omega^1_2 \\
\hline
\omega^2_1 \\
\omega^2_2 \\
\hline
\omega^3_1 \\
\omega^3_2 \\
\vdots
\end{array}\right)
\end{eqnarray*}
Now the first five rows match $\bfS$. One more operation
\[
R_6 \rightarrow \D\frac{R_6 - (\dot{u}_2R_4 + u_1R_1 + x_1R_5)}{u_2}
\]
puts the coframing in the following form
\begin{eqnarray*}
& \left(\begin{array}{c}
\Omega^0_2 \\
\Omega^0_3 \\
x_1\Omega^0_3 - \Omega^0_1 \\
\hline
\left( \Omega^1_1 - x_1\Omega^0_2 \right) / u_2 \\
\Omega^1_2 \\
\hline
\Big( \Omega^2_1 - \left( \dot{u}_2\left( u_2\Omega^1_1 + x_1\Omega^0_2 \right) + u_1\Omega^0_2 + x_1\Omega^1_2 \right) \Big) / u_2 \\
\Omega^2_2 \\
\vdots
\end{array}\right) \\
& = \left(\begin{array}{ccc|cc|cc|ccc}
0 & 0 & 0 & 0 & 1 & 0 & 0 & 0 & 0 & \cdots \\
0 & 1 & 0 & 0 & 0 & 0 & 0 & 0 & 0 & \cdots \\
0 & 0 & 1 & 0 & 0 & 0 & 0 & 0 & 0 & \cdots \\
\hline
1 & 0 & 0 & 0 & 0 & 0 & 0 & 0 & 0 & \cdots \\
0 & 0 & 0 & 0 & 0 & 0 & 1 & 0 & 0 & \cdots \\
\hline
0 & 0 & 0 & 1 & 0 & 0 & 0 & 0 & 0 & \cdots \\
0 & 0 & 0 & 0 & 0 & 0 & 0 & 0 & 1 & \cdots \\
\vdots & \vdots & \vdots & \vdots & \vdots & \vdots & \vdots & \vdots & \vdots \\
\end{array}\right)
\left(\begin{array}{c}
\omega^0_1 \\
\omega^0_2 \\
\omega^0_3 \\
\hline
\omega^1_1 \\
\omega^1_2 \\
\hline
\omega^2_1 \\
\omega^2_2 \\
\hline
\omega^3_1 \\
\omega^3_2 \\
\vdots
\end{array}\right),
\end{eqnarray*}
and all visible rows now match those of $\bfS$. Continuing this process ad infinitum gives us new coframings that transform via $\bfS$. At present, this transformation looks like
\[
\bfgg^{-1}\bfO = \bfS \bfG \bfo,
\]
where $\bfG$ is the identity $Id$. To put it in the desired form, we simply invert the action on the left hand side. This results in the following factored transformation.
\begin{eqnarray*}
\bfO = \left(\begin{array}{ccc|cc|ccc}
0 & x_1 & -1 & 0 & 0 & 0 & 0 & \cdots \\
1 & 0 & 0 & 0 & 0 & 0 & 0 & \cdots \\
0 & 1 & 0 & 0 & 0 & 0 & 0 & \cdots \\
\hline
x_1 & 0 & 0 & u_2 & 0 & 0 & 0 & \cdots \\
0 & 0 & 0 & 0 & 1 & 0 & 0 & \cdots \\
\hline
-x_1\dot{u}_2 - u_1 & 0 & 0 & u_2\dot{u}_2 & x_1 & u_2 & 0 & \cdots \\
0 & 0 & 0 & 0 & 0 & 0 & 1 & \cdots \\
\vdots & \vdots & \vdots & \vdots & \vdots & \vdots & \vdots \\
\end{array}\right)
\bfS\ Id\ \bfo
\end{eqnarray*}
Note that this decomposition using non-autonomous group elements is not unique, however it was chosen so that the second non-autonomous group element of the latter equation was particularly simple (the identity in this case). Any problem with three states and two controls can be simplified in a similar way, as we will see below.

\section{Three States and Two Controls}
\label{3states}

\subsection{Preferred Structure Equations}

In the method of equivalence, described in Chapter \ref{equivalence}, one important step is to work with an initial, preferred coframing that encapsulates the problem at hand and satisfies some particularly nice relations that ought to be preserved by the equivalence in question. In this section we will make one final refinement to our coframings \eqref{coframe} so that they satisfy some particularly nice structure equations that ought to be preserved by dynamic equivalence.

Note that for a control system $\dot{\bfx} = \bff(\bfx,\bfu)$ with $n$ state variables and $s\leq n$ control variables, the vector
\begin{eqnarray*}
\bff(\bfx,\bfu) & = & \left(\begin{array}{c}
f_1(\bfx,\bfu) \\
\vdots \\
f_n(\bfx,\bfu)
\end{array}\right)
\end{eqnarray*}
must have $\textrm{rank}\frac{\partial\bff}{\partial\bfu} = s$. Therefore, by the implicit function theorem, a static equivalence always exists so that the above system is equivalent to $\dot{\tilde\bfx}=\tilde\bff(\tilde\bfx,\tilde\bfu)$ where
\begin{eqnarray*}
\tilde\bff(\tilde\bfx,\tilde\bfu) & = & \left(\begin{array}{c}
\tilde{u}_1 \\
\vdots \\
\tilde{u}_s \\
\tilde{f}_{s+1}(\tilde\bfx,\tilde\bfu) \\
\vdots \\
\tilde{f}_n(\tilde\bfx,\tilde\bfu)
\end{array}\right),
\end{eqnarray*}
where $\tilde{x}_i = x_i$ for $1\leq i\leq n$ up to reordering and $\tilde{u}_j = f_j(\bfx,\bfu)$ for $1\leq j\leq  s$.

We will now, and for the rest of the paper, concern ourselves with the case of three state variables and two control variables. The above adaptation suggests altering \eqref{coframe} for the case of three states and two controls to the following coframing.

\begin{eqnarray*}
\tilo^{-1} & = & dt \\
\tilo^0_1 & = & dx_1 - u_1\ dt\\
\tilo^0_2 & = & dx_2 - u_2\ dt\\
\tilo^0_3 & = & dx_3 - f(\bfx,\bfu)\ dt \\
\tilo^1_1 & = & du_1 - \dot{u}_1\ dt\\
\tilo^1_2 & = & du_2 - \dot{u}_2\ dt\\
\vdots
\end{eqnarray*}

Here $f(\bfx,\bfu)$ is a scalar function. Note that in this coframing, $d\tilo^i_j = -\tilo^{i+1}_j\wedge\tilo^{-1}$ for $i\geq 0$ and $j=1,2$. The outlier in this nice pattern of exterior derivatives is, of course,
\[
d\tilo^0_3 = - \sum_{i=1}^3 f_{x_i}(x,u) \tilo^0_i\wedge\tilo^{-1} - \sum_{i=1}^2 f_{u_i}(x,u) \tilo^1_i\wedge\tilo^{-1}.
\]
With one more adaptation of the coframing, we can make even this structure equation easier to work with. Let the following be our preferred coframing for the case of $n=3$ state variables, $s=2$ control variables.
\begin{eqnarray}
\label{final_coframe}
\begin{array}{rl}
\omega^{-1} & =\ dt \\
\omega^0_1 & =\ dx_1 - u_1\ dt\\
\omega^0_2 & =\ dx_2 - u_2\ dt\\
\omega^0_3 & =\ dx_3 - f\ dt - f_{u_1} (dx_1 - u_1\ dt) - f_{u_2}(dx_2 - u_2\ dt) \\
\omega^1_1 & =\ du_1 - \dot{u}_1\ dt\\
\omega^1_2 & =\ du_2 - \dot{u}_2\ dt\\
\vdots
\end{array}
\end{eqnarray}
Note that this coframing satisfies some particularly nice structure equations.
\begin{eqnarray*}
\begin{array}{rl}
d{\omega}^0_1 & = -\omega^1_1 \wedge \omega^{-1} \\
d{\omega}^0_2 & = -\omega^1_2 \wedge \omega^{-1} \\
d{\omega}^0_3 & \equiv 0 \qquad \textrm{ mod }\ \bfo^0 \\
d{\omega}^j_k & = -\omega^{j+1}_k \wedge \omega^{-1} \qquad (j>0,\ k=1,2)
\end{array}
\end{eqnarray*}
We will take this coframing, along with the analogous coframing $\bfO$ in $(\bfy,\bfv)$ coordinates, as our starting point. Let $\bfbo = \bfG\bfo$ and $\bfbO = \bfS\bfbo$ so that $\bfO = \bfgg\bfbO$. In addition, we will require that at every step of our transformation of the coframes, $\bfbo$, $\bfbO$ preserves the following nice properties of the structure equations and their algebraic ideals:
\begin{eqnarray}
\label{struc_eqns}
\begin{array}{l}
\left. \begin{array}{rl}
d{\omega}^0_1 & \equiv -\omega^1_1 \wedge \omega^{-1} \\
d{\omega}^0_2 & \equiv -\omega^1_2 \wedge \omega^{-1} \\
d{\omega}^0_3 & \equiv 0
\end{array} \right\} \quad \textrm{ mod }\ \bfo^0 \\
d{\omega}^j_k \equiv -\omega^{j+1}_k \wedge \omega^{-1} \quad \textrm{ mod }\ \{\ \bfo^i\ |\ 0\leq i\leq j\ \}, \quad (j>0,\ k=1,2).
\end{array}
\end{eqnarray}

\subsection{Reducing $\bfG$}

Consider the coframing $\bfbO = \bfS\bfG\bfo$. Since we plan on applying a generic $\bfgg$ in the non-autonomous problem $\bfO = \bfgg\bfbO$, $\bfG$ does not have to be completely generic. It can be simplified to remove some redundancies. For example, since $\bo^0_3 \mapsto \bO^0_3$ under $\bfS$, there is no need to add an arbitrary multiple of $\bo^0_3$ to any other form through $\bfG$ since this can be taken care of with $\bfgg$. What follows will illustrate this more explicitly.

We have coframings $\bfbo = \bfG\bfo$ and $\bfbO = \bfS\bfbo = \bfS\bfG\bfo$. Recall that $G^i_i = G^1_1$ for all $i\geq 1$ by Theorem \ref{Arepeats}. Consider the following identities.
\begin{eqnarray*}
\bO^0_1 & = & (G^1_0)^2_1\ \omega^0_1 + (G^1_0)^2_2\ \omega^0_2 + (G^1_0)^2_3\ \omega^0_3 + (G^1_1)^2_1\ \omega^1_1 + (G^1_1)^2_2\ \omega^1_2 \\
\bO^0_2 & = & (G^0_0)^2_1\ \omega^0_1 + (G^0_0)^2_2\ \omega^0_2 + (G^0_0)^2_3\ \omega^0_3 \\
\bO^0_3 & = & (G^0_0)^3_1\ \omega^0_1 + (G^0_0)^3_2\ \omega^0_2 + (G^0_0)^3_3\ \omega^0_3 \\
\\
\bO^1_1 & = & (G^0_0)^1_1\ \omega^0_1 + (G^0_0)^1_2\ \omega^0_2 + (G^0_0)^1_3\ \omega^0_3 \\
\bO^1_2 & = & (G^2_0)^2_1\ \omega^0_1 + (G^2_0)^2_2\ \omega^0_2 + (G^2_0)^2_3\ \omega^0_3 + (G^2_1)^2_1\ \omega^1_1 + (G^2_1)^2_2\ \omega^1_2 \\
& & +\ (G^1_1)^2_1\ \omega^2_1 + (G^1_1)^2_2\ \omega^2_2 \\
\\
\bO^2_1 & = & (G^1_0)^1_1\ \omega^0_1 + (G^1_0)^1_2\ \omega^0_2 + (G^1_0)^1_3\ \omega^0_3 + (G^1_1)^1_1\ \omega^1_1 + (G^1_1)^1_2\ \omega^1_2 \\
\bO^2_2 & = & (G^3_0)^2_1\ \omega^0_1 + (G^3_0)^2_2\ \omega^0_2 + (G^3_0)^2_3\ \omega^0_3 + (G^3_1)^2_1\ \omega^1_1 + (G^3_1)^2_2\ \omega^1_2 \\
& & +\ (G^3_2)^2_1\ \omega^2_1 + (G^3_2)^2_2\ \omega^2_2 + (G^1_1)^2_1\ \omega^3_1 + (G^1_1)^2_2\ \omega^3_2 \\
\\
\bO^3_1 & = & (G^2_0)^1_1\ \omega^0_1 + (G^2_0)^1_2\ \omega^0_2 + (G^2_0)^1_3\ \omega^0_3 + (G^2_1)^1_1\ \omega^1_1 + (G^2_1)^1_2\ \omega^1_2 \\
& & +\ (G^1_1)^1_1\ \omega^2_1 + (G^1_1)^1_2\ \omega^2_2 \\
\bO^3_2 & = & (G^4_0)^2_1\ \omega^0_1 + (G^4_0)^2_2\ \omega^0_2 + (G^4_0)^2_3\ \omega^0_3 + (G^4_1)^2_1\ \omega^1_1 + (G^4_1)^2_2\ \omega^1_2 \\
& & +\ (G^4_2)^2_1\ \omega^2_1 + (G^4_2)^2_2\ \omega^2_2 + (G^4_3)^2_1\ \omega^3_1 + (G^4_3)^2_2\ \omega^3_2 + (G^1_1)^2_1\ \omega^4_1 + (G^1_1)^2_2\ \omega^4_2 \\
 & \vdots
\end{eqnarray*}
Now $\bfgg$ will add arbitrary multiples of $\bO^0_2$ and $\bO^0_3$ to every other part of the coframing in order to get the final coframing $\bfO$. Since they are linearly independent, they do not need to be completely arbitrary. We will not lose anything by letting $(G^0_0)^2_2 = (G^0_0)^3_3 = 1$ and $(G^0_0)^3_2 = (G^0_0)^2_3 = 0$. In fact all of the other terms above involving $\omega^0_2$ and $\omega^0_3$ may as well be set to zero since $\bfgg$ will take care of these through nonautonomous static equivalence.
\begin{eqnarray*}
\bO^0_1 & = & (G^1_0)^2_1\ \omega^0_1 + (G^1_1)^2_1\ \omega^1_1 + (G^1_1)^2_2\ \omega^1_2 \\
\bO^0_2 & = & (G^0_0)^2_1\ \omega^0_1 + \omega^0_2 \\
\bO^0_3 & = & (G^0_0)^3_1\ \omega^0_1 + \omega^0_3 \\
\\
\bO^1_1 & = & (G^0_0)^1_1\ \omega^0_1 \\
\bO^1_2 & = & (G^2_0)^2_1\ \omega^0_1 + (G^2_1)^2_1\ \omega^1_1 + (G^2_1)^2_2\ \omega^1_2 + (G^1_1)^2_1\ \omega^2_1 + (G^1_1)^2_2\ \omega^2_2 \\
\\
\bO^2_1 & = & (G^1_0)^1_1\ \omega^0_1 + (G^1_1)^1_1\ \omega^1_1 + (G^1_1)^1_2\ \omega^1_2 \\
\bO^2_2 & = & (G^3_0)^2_1\ \omega^0_1 + (G^3_1)^2_1\ \omega^1_1 + (G^3_1)^2_2\ \omega^1_2 + (G^3_2)^2_1\ \omega^2_1 + (G^3_2)^2_2\ \omega^2_2 \\
& & +\ (G^1_1)^2_1\ \omega^3_1 + (G^1_1)^2_2\ \omega^3_2 \\
\\
\bO^3_1 & = & (G^2_0)^1_1\ \omega^0_1 + (G^2_1)^1_1\ \omega^1_1 + (G^2_1)^1_2\ \omega^1_2 + (G^1_1)^1_1\ \omega^2_1 + (G^1_1)^1_2\ \omega^2_2 \\
\bO^3_2 & = & (G^4_0)^2_1\ \omega^0_1 + (G^4_1)^2_1\ \omega^1_1 + (G^4_1)^2_2\ \omega^1_2 + (G^4_2)^2_1\ \omega^2_1 + (G^4_2)^2_2\ \omega^2_2 \\
& & +\ (G^4_3)^2_1\ \omega^3_1 + (G^4_3)^2_2\ \omega^3_2 +\ (G^1_1)^2_1\ \omega^4_1 + (G^1_1)^2_2\ \omega^4_2 \\
 & \vdots
\end{eqnarray*}
Of course we are keeping careful note that every group reduction we have made is allowed due to the freedom we have in choosing $\bfgg$.

Now it is clear that we may as well choose $(G^0_0)^1_1 = 1$, and thus we may also set any term involving $\omega^0_1$ below $\bO^1_1$ to zero since $\bfgg$ will be adding arbitrary multiples of $\bO^1_1$ to these.
\begin{eqnarray*}
\bO^0_1 & = & (G^1_0)^2_1\ \omega^0_1 + (G^1_1)^2_1\ \omega^1_1 + (G^1_1)^2_2\ \omega^1_2 \\
\bO^0_2 & = & (G^0_0)^2_1\ \omega^0_1 + \omega^0_2 \\
\bO^0_3 & = & (G^0_0)^3_1\ \omega^0_1 + \omega^0_3 \\
\\
\bO^1_1 & = & \omega^0_1 \\
\bO^1_2 & = & (G^2_1)^2_1\ \omega^1_1 + (G^2_1)^2_2\ \omega^1_2 + (G^1_1)^2_1\ \omega^2_1 + (G^1_1)^2_2\ \omega^2_2 \\
\\
\bO^2_1 & = & (G^1_1)^1_1\ \omega^1_1 + (G^1_1)^1_2\ \omega^1_2 \\
\bO^2_2 & = & (G^3_1)^2_1\ \omega^1_1 + (G^3_1)^2_2\ \omega^1_2 + (G^3_2)^2_1\ \omega^2_1 + (G^3_2)^2_2\ \omega^2_2 + (G^1_1)^2_1\ \omega^3_1 + (G^1_1)^2_2\ \omega^3_2 \\
\\
\bO^3_1 & = & (G^2_1)^1_1\ \omega^1_1 + (G^2_1)^1_2\ \omega^1_2 + (G^1_1)^1_1\ \omega^2_1 + (G^1_1)^1_2\ \omega^2_2 \\
\bO^3_2 & = & (G^4_1)^2_1\ \omega^1_1 + (G^4_1)^2_2\ \omega^1_2 + (G^4_2)^2_1\ \omega^2_1 + (G^4_2)^2_2\ \omega^2_2 + (G^4_3)^2_1\ \omega^3_1 + (G^4_3)^2_2\ \omega^3_2 \\
 & & +\ (G^1_1)^2_1\ \omega^4_1 + (G^1_1)^2_2\ \omega^4_2 \\
 & \vdots
\end{eqnarray*}
One entry in every $\bO^i_j$ can be scaled to unity. Note that $G^1_1$ is an invertible $2\times 2$ matrix, so that either the pair $(G^1_1)^1_1$, $(G^1_1)^2_2$ or $(G^1_1)^1_2$, $(G^1_1)^2_1$ is nonzero. If the former pair is zero, then $\bfgg$ would allow us to switch the roles of every $\bO^i_1$ and $\bO^i_2$ for $i\geq 1$. Thus without loss of generality we can let $(G^1_1)^1_1 = (G^1_1)^2_2 = 1$. The arbitrariness of $\bfgg$ will then let us cancel out any terms below these scaled terms. For example, adding multiples of $\bO^0_1$ and $\bO^1_1$ to $\bO^1_2$ will get rid of the $\omega^1_2$ term in all the $\bO^i$, $i\geq 1$. We can also scale the $\omega^2_2$ term in $\bO^1_2$ to unity, and thus every $\omega^2_2$ below can be eliminated. After this process of scaling one term per $\bO^i_j$ and using this to eliminate the appropriate terms below, we are left with the following.
\begin{eqnarray*}
\bO^0_1 & = & (G^1_0)^2_1\ \omega^0_1 + (G^1_1)^2_1\ \omega^1_1 + \omega^1_2 \\
\bO^0_2 & = & (G^0_0)^2_1\ \omega^0_1 + \omega^0_2 \\
\bO^0_3 & = & (G^0_0)^3_1\ \omega^0_1 + \omega^0_3 \\
\\
\bO^1_1 & = & \omega^0_1 \\
\bO^1_2 & = & (G^2_1)^2_1\ \omega^1_1 + (G^1_1)^2_1\ \omega^2_1 + \omega^2_2 \\
\\
\bO^2_1 & = & \omega^1_1 \\
\bO^2_2 & = & (G^3_2)^2_1\ \omega^2_1 + (G^1_1)^2_1\ \omega^3_1 + \omega^3_2 \\
\\
\bO^3_1 & = & \omega^2_1 \\
\bO^3_2 & = & (G^4_3)^2_1\ \omega^3_1 + (G^1_1)^2_1\ \omega^4_1 + \omega^4_2 \\
 & \vdots
\end{eqnarray*}

After all such redundancies are removed, this is what our group element, now called $\bfGG$, looks like.
\[
\bfGG = 
\left( \begin{array}{c|ccc|cc|cc|cc|ccc}
1 & 0 & 0 & 0 & 0 & 0 & 0 & 0 & 0 & 0 & 0 & 0 & \cdots  \\
\hline
0 & 1 & 0 & 0 & 0 & 0 & 0 & 0 & 0 & 0 & 0 & 0 & \cdots  \\
0 & (G^0_0)^2_1 & 1 & 0 & 0 & 0 & 0 & 0 & 0 & 0 & 0 & 0 & \cdots  \\
0 & (G^0_0)^3_1 & 0 & 1 & 0 & 0 & 0 & 0 & 0 & 0 & 0 & 0 & \cdots  \\
\hline
0 & 0 & 0 & 0 & 1 & 0 & 0 & 0 & 0 & 0 & 0 & 0 & \cdots  \\
0 & (G^1_0)^2_1 & 0 & 0 & (G^1_1)^2_1 & 1 & 0 & 0 & 0 & 0 & 0 & 0 & \cdots  \\
\hline
0 & 0 & 0 & 0 & 0 & 0 & 1 & 0 & 0 & 0 & 0 & 0 & \cdots  \\
0 & 0 & 0 & 0 & (G^2_1)^2_1 & 0 & (G^1_1)^2_1 & 1 & 0 & 0 & 0 & 0 & \cdots  \\
\hline
0 & 0 & 0 & 0 & 0 & 0 & 0 & 0 & 1 & 0 & 0 & 0 & \cdots  \\
0 & 0 & 0 & 0 & 0 & 0 & (G^3_2)^2_1 & 0 & (G^1_1)^2_1 & 1 & 0 & 0 & \cdots  \\
 & & & & \vdots \\
\end{array} \right)
\]

We will need to use the fact that $\bfbo$ is a coframing. Therefore exterior derivatives of the entries of $\bfGG$ can be written as linear combinations of these. Note that as far as we know, every $d(G^i_j)^k_l$ could be linear combinations of $\bfbo^r$ for some unknown $r$. We will employ the following notation:
\[
d(G^i_j)^k_l = (G^i_j)^{k}_{l,-1} \omega^{-1} + \sum_\alpha\sum_\beta (G^i_j)^{k,\beta}_{l,\alpha} \omega^\alpha_\beta
\]
We will show below that $r$ is not arbitrarily large by looking at structure equations.

By investigating $d\bfbO$, we can further reduce the entries of $\bfGG$. Until stated otherwise, the following equivalences $\equiv$ are modulo $\bO^0_i$, $i=1,2,3$. We will start with $d\bO^0_3$.
\begin{eqnarray*}
\bO^0_3 & = & \omega^0_3 + (G^0_0)^3_1\omega^0_1 \\
d\bO^0_3 & = & d\omega^0_3 + d(G^0_0)^3_1\wedge\omega^0_1 + (G^0_0)^3_1\ d\omega^0_1 \\
 & \equiv & \Big[ \Big( u_2 f_{u_1x_2} - f_{x_1} - (G^0_0)^3_{1,-1} + u_1 f_{u_1x_1} + f_{u_1x_3} f + f_{x_3} (G^0_0)^3_1 - f_{x_3} f_{u_1} \\
& & +\ \dot{u}_1 f_{u_1u_1} + \dot{u}_2 f_{u_1u_2} \Big) - (G^0_0)^2_1\Big( u_2 f_{u_2x_2} + u_1 f_{u_2x_1} + \dot{u}_2 f_{u_2u_2} - f_{x_2} \\
& & +\ f_{u_2x_3} f + \dot{u}_1 f_{u_1u_2} - f_{x_3} f_{u_2} \Big) \Big]\ \bO^1_1\wedge\bO^{-1} + \Big[ \Big( f_{u_1u_2} - f_{u_2u_2} (G^1_1)^2_1 \Big)(G^0_0)^2_1 \\
& & -\ \Big( (G^0_0)^{3,1}_{1,2} (G^1_1)^2_1 - f_{u_1u_2} (G^1_1)^2_1 + f_{u_1u_1} - (G^0_0)^{3,1}_{1,1} \Big) \Big] \bO^2_1\wedge\bO^1_1 \\
& & -\ (G^0_0)^3_1 \ \bO^2_1\wedge\bO^{-1}
\end{eqnarray*}
Anything above that is not a multiple of $\bO^1_1\wedge\bO^{-1}$ or $\bO^1_2\wedge\bO^{-1}$ must have zero coefficient. Of greatest interest at the moment is the term $\bO^2_1\wedge\bO^{-1}$. Since this cannot be here, its coefficient must be zero.
\begin{eqnarray}
\label{G00_31}
(G^0_0)^3_1 = 0
\end{eqnarray}
There is also a $\bO^2_1\wedge\bO^1_1$ term which must vanish. Through the above equation, this simplifies to the following.
\[
\Big( f_{u_1u_2} - f_{u_2u_2} (G^1_1)^2_1 \Big)(G^0_0)^2_1 + \Big( f_{u_1u_2} (G^1_1)^2_1 - f_{u_1u_1} \Big) = 0
\]

Moving on, we will look at $d\bO^0_2$.
\begin{eqnarray*}
\bO^0_2 & = & \omega^0_2 + (G^0_0)^2_1\ \omega^0_1 \\
d\bO^0_2 & = & d\omega^0_2 + d(G^0_0)^2_1\wedge\omega^0_1 + (G^0_0)^2_1\ d\omega^0_1\\
& \equiv & \Big( (G^1_0)^2_1 - (G^0_0)^2_{1,-1} \Big)\ \bO^1_1\wedge\bO^{-1} + \Big( (G^1_1)^2_1 - (G^0_0)^2_1 \Big)\  \bO^2_1\wedge\bO^{-1} \\
 & & \qquad +\ \sum_{i=1}^\infty \Big( (G^0_0)^{2,1}_{1,i} - (G^0_0)^{2,2}_{1,i} (G^1_1)^2_1 - (G^0_0)^{2,2}_{1,i+1} (G^{i+1}_i)^2_1 \Big)\ \bO^{i+1}_1\wedge\bO^1_1 \\
 & & \qquad \qquad +\ \sum_{i=1}^\infty (G^0_0)^{2,2}_{1,i+1}\bO^i_2\wedge\bO^1_1
\end{eqnarray*}
Similarly here it is the vanishing of the $\bO^2_1\wedge\bO^{-1}$ term that tells us
\[
(G^1_1)^2_1 = (G^0_0)^2_1.
\]
The vanishing of the terms in the final two summations tells us
\begin{eqnarray*}
(G^0_0)^{2,1}_{1,1} & = & (G^0_0)^{2,1}_{1,2} (G^0_0)^2_1, \\
(G^0_0)^{2,1}_{1,i} & = & 0, \\
(G^0_0)^{2,2}_{1,i} & = & 0
\end{eqnarray*}
for all $i\geq 2$. We knew that
\[
d(G^0_0)^2_1 = (G^0_0)^{2}_{1,-1} \omega^{-1} + \sum_\alpha\sum_\beta (G^0_0)^{2,\beta}_{1,\alpha} \omega^\alpha_\beta
\]
had to be a finite sum, and now we have a bound on where that sum must terminate ($\alpha=1$).

Now consider $d\bO^0_1$.
\begin{eqnarray*}
\bO^0_1 & = & (G^1_0)^2_1\ \omega^0_1 + (G^0_0)^2_1\ \omega^1_1 + \omega^1_2 \\
d\bO^0_1 & = & d(G^1_0)^2_1\wedge \omega^0_1 + (G^1_0)^2_1\ d\omega^0_1 + d(G^0_0)^2_1\wedge\omega^1_1 + (G^0_0)^2_1\ d\omega^1_1 + d\omega^1_2 \\
& \equiv & \Big[ \Big( (G^0_0)^{2,2}_{1,0} - (G^1_0)^{2,2}_{1,1} \Big) (G^0_0)^2_1 + (G^1_0)^{2,1}_{1,1} - (G^0_0)^{2,1}_{1,0} - (G^1_0)^{2,2}_{1,2} (G^2_1)^2_1 \\
& & +\ (G^0_0)^{2,2}_{1,1} (G^1_0)^2_1 \Big]\ \bO^2_1\wedge\bO^1_1 - (G^1_0)^2_{1,-1} \bO^1_1\wedge\bO^{-1} - \bO^1_2\wedge\bO^{-1} \\
& & -\ (G^1_0)^{2,2}_{1,2} \bO^1_1\wedge\bO^1_2 + \Big( (G^2_1)^2_1 - (G^1_0)^2_1 - (G^0_0)^2_{1,-1} \Big)\ \bO^2_1\wedge\bO^{-1} \\
& & +\ \sum_{i=2}^\infty \Big( (G^1_0)^{2,1}_{1,i} - (G^1_0)^{2,2}_{1,i} (G^0_0)^2_1 - (G^1_0)^{2,2}_{1,i+1} (G^{i+1}_i)^2_1 \Big)\ \bO^{i+1}_1\wedge\bO^1_1 \\
& & +\ \sum_{i=2}^\infty (G^1_0)^{2,2}_{1,i+1}\bO^i_2\wedge\bO^1_1
\end{eqnarray*}
The relations that come from this calculation are these for $i\geq 2$.
\begin{eqnarray*}
(G^2_1)^2_1 & = & (G^1_0)^2_1 + (G^0_0)^2_{1,-1} \\
(G^1_0)^{2,1}_{1,1} & = &  (G^1_0)^{2,2}_{1,1} (G^0_0)^2_1 + (G^0_0)^{2,1}_{1,0} - (G^0_0)^{2,2}_{1,0} (G^0_0)^2_1 - (G^0_0)^{2,2}_{1,1} (G^1_0)^2_1 \\
(G^1_0)^{2,2}_{1,i} & = & 0 \\
(G^1_0)^{2,1}_{1,i} & = & 0 
\end{eqnarray*}
Therefore we have found a bound on the sum for $d(G^1_0)^{2,2}$ as well.

Continuing this process for higher order terms yields the following important result
\[
(G^{i+2}_{i+1})^2_1 = (G^1_0)^2_1 + (i+1)(G^0_0)^2_{1,-1}
\]
for $i\geq 1$.

To review, $\bfGG$ now has the following form.
\begin{multline}
\label{Gnice}
\bfGG= \\
\left( \begin{array}{c|c@{\hspace{1pt}}c@{\hspace{1pt}}c|c@{\hspace{1pt}}c|c@{\hspace{1pt}}c|c@{\hspace{3pt}}c|c@{\hspace{3pt}}c@{\hspace{3pt}}c}
1 & 0 & 0 & 0 & 0 & 0 & 0 & 0 & 0 & 0 & 0 & 0 & \cdots  \\
\hline
0 & 1 & 0 & 0 & 0 & 0 & 0 & 0 & 0 & 0 & 0 & 0 & \cdots  \\
0 & (G^0_0)^2_1 & 1 & 0 & 0 & 0 & 0 & 0 & 0 & 0 & 0 & 0 & \cdots  \\
0 & 0 & 0 & 1 & 0 & 0 & 0 & 0 & 0 & 0 & 0 & 0 & \cdots  \\
\hline
0 & 0 & 0 & 0 & 1 & 0 & 0 & 0 & 0 & 0 & 0 & 0 & \cdots  \\
0 & (G^1_0)^2_1 & 0 & 0 & (G^0_0)^2_1 & 1 & 0 & 0 & 0 & 0 & 0 & 0 & \cdots  \\
\hline
0 & 0 & 0 & 0 & 0 & 0 & 1 & 0 & 0 & 0 & 0 & 0 & \cdots  \\
0 & 0 & 0 & 0 & (G^1_0)^2_1\! +\! (G^0_0)^2_{1,-1} & 0 & (G^0_0)^2_1 & 1 & 0 & 0 & 0 & 0 & \cdots  \\
\hline
0 & 0 & 0 & 0 & 0 & 0 & 0 & 0 & 1 & 0 & 0 & 0 & \cdots  \\
0 & 0 & 0 & 0 & 0 & 0 & (G^1_0)^2_1\! +\! 2(G^0_0)^2_{1,-1} & 0 & (G^0_0)^2_1 & 1 & 0 & 0 & \cdots  \\
 & & & & \vdots & & & & & & \\
\end{array} \right)
\end{multline}

What we have boiled the problem down to now is the equivalence $\bfO = \bfgg\bfbO$, where the coframing $\bfbO$ contains three functions $f$, $(G^0_0)^2_1$, and $(G^1_0)^2_1$.

\underline{Remark:} An important but subtle point to take note of is the following: we have singled out $\bO^0_1$ through $\bfS$ as the piece of the coframing $\bfbO^0$ that contains higher order terms in $\bfo$, and we have also singled out $\bO^0_3$ by choosing an adapted coframing with $d\bO^0_3 \equiv 0$ mod $\bfbO^0$, and these two choices are compatible.

This fact is actually quite easy to see. In our coframings, note that $\bO^0_3 = \omega^0_3$. Since $\bfgg$ preserves the span of $\{\ \bO^0_1,\bO^0_2,\bO^0_3\ \}$, $\omega^0_3$ must be in the span of $\{\ \Omega^0_1,\Omega^0_2,\Omega^0_3\ \}$. Thus $\omega^0_3$, which has the property that $d\omega^0_3\equiv 0$ mod $\bfo^0$, does not also get bumped up in the dynamically equivalent coframing to a higher order term.

\section{Dynamic Equivalence of Control Affine Systems}
\label{affinesystems}

Keep in mind at this point that we are concerned with dynamic equivalence, which is a weaker equivalence than static equivalence. The static equivalence case was dealt with first in the control linear case of three states and two controls by Wilkens and later by Elkin in the affine linear case up to four states. The representatives of the five distinct static equivalent control affine systems with three states and two controls put forth by Elkin are these:
\begin{eqnarray}
\begin{array}{rcl}
\dot{x}_1 & = & u_1 \\
\dot{x}_2 & = & u_2 \\
\dot{x}_3 & = & f(\bfx,\bfu)
\end{array}
\end{eqnarray}
where $f(\bfx,\bfu)$ is one of the five following functions:
\begin{eqnarray*}
\label{affine}
\left\{ \begin{array}{c}
0 \\
1 \\
x_2 \\
x_2 u_1 \\
1 + x_2 u_1
\end{array}\right\}
\end{eqnarray*}

In this section, we will finally put to use our previous results involving infinite prolongations and the factorization of coframing pullbacks. We show, using arguments about certain ideals preserved under dynamic equivalence, that neither of the first two systems listed above are dynamically equivalent to any other control system with $J=K=0$. The proof of the final theorem gives explicit dynamic equivalences between the last three systems above.

\begin{theorem}
The control system corresponding to $\dot{x}_3=0$ with two control variables is not dynamically equivalent to any other control system with $J=K=0$ to which it is not static equivalent.
\end{theorem}

\begin{proof}
Suppose $\bfO = \bfgg\bfS\bfGG\bfo$, where $\bfGG$ is given by \eqref{Gnice}, $\bfS$ is given by \eqref{stackpole}, $\bfgg$ is a generic nonautonomous static equivalence, and $\bfo$ is the following coframing for $\dot{x}_3=0$.
\begin{eqnarray*}
\omega^{-1} & = & dt \\
\omega^0_1 & = & dx_1 - u_1\ dt \\
\omega^0_2 & = & dx_2 - u_2\ dt \\
\omega^0_3 & = & dx_3 \\
& \vdots &
\end{eqnarray*}
The coframing $\bfbO = \bfS\bfGG\bfo$ would then look as follows.
\begin{eqnarray*}
\bO^{-1} & = & dt \\
\bO^0_1 & = & \left( (G^1_0)^2_1\ dx_1 + (G^0_0)^2_1\ du_1 + du_2 \right) - \left( (G^1_0)^2_1\ u_1 + (G^0_0)^2_1\ \dot{u}_1 + \dot{u}_2 \right)\ dt \\
\bO^0_2 & = & \left( dx_2 + (G^0_0)^2_1\ dx_1 \right) - \left( (G^0_0)^2_1\ u_1 + u_2 \right)\ dt \\
\bO^0_3 & = & dx_3 \\
& \vdots &
\end{eqnarray*}

Now notice that the algebraic ideal $\bfbO^0$ is preserved by $\bfgg$. But all of our equivalences also preserve $t$, and hence $dt$. Therefore, if $\overline{\Lambda}_\infty$ has the coframing
\begin{eqnarray*}
\Omega^{-1} & = & dt \\
\Omega^0_1 & = & dy_1 - g_1(y,v)\ dt \\
\Omega^0_2 & = & dy_2 - g_2(y,v)\ dt \\
\Omega^0_3 & = & dy_3 - g_3(y,v)\ dt \\
& \vdots &
\end{eqnarray*}
we get that $\{\bO^0_1,\bO^0_2,\bO^0_3\} \equiv \{dy_1,dy_2,dy_3\}$ mod $dt$. Since this is an integrable ideal that contains $\bO^0_3=dx_3$, we can arrange through the appropriate choice of $\bfgg$ that $dy_3=dx_3$. Note that this automatically satisfies $d\bO^0_3 \equiv 0$ mod $\bO^0_1,\bO^0_2,\bO^0_3$ since $d\bO^0_3$ is identically zero.

Therefore $\dot{y}_3=\dot{x}_3=0$. What we have done is taken any strict dynamic equivalence to the system $\dot{x}_3=0$ with $J=K=0$ and altered it via static equivalence to a strict dynamic equivalence to itself. So any control system that is dynamically equivalent to $\dot{x}_3=0$ with $J=K=0$ is in fact a dynamic equivalence to a system that is static equivalent to $\dot{x}_3=0$.
\end{proof}

\begin{theorem}
The control system corresponding to $\dot{x}_3=1$ with two control variables is not dynamically equivalent to any other control system with $J=K=0$ to which it is not static equivalent.
\end{theorem}

\begin{proof}
The proof is nearly identical to that of the previous theorem. Replace $\dot{x}_3=0$ with $\dot{x}_3=1$, and proceed in the same fashion.
\end{proof}

Note that the method used in the previous two theorems could also be applied to the case of $\dot{x}_3=x_2$. A difference occurs, however, when reaching the step $\dot{y}_3=\dot{x}_3=x_2$. Since $x_2$ is not necessarily equal to $y_2$, we see that the resulting system may or may not necessarily be static equivalent to the original system $\dot{x}_3=x_2$. It in fact turns out, as stated in the next theorem, that this new system need not be static equivalent to the original system.

\begin{theorem}
The control systems $\dot{x}_3 = x_2$, $x_2u_1$, $1+x_2u_1$ are strictly dynamically equivalent to each other.
\end{theorem}

\begin{proof}
The following sets of maps between infinite jet bundles give explicit dynamic equivalences for the three systems. We will demonstrate that the maps take solutions of one control system to solutions of the other. The fact that the maps composed with their respective inverses are in fact the identity on solutions is simple enough and is left to the reader. 
\[
\begin{array}{ccccc}
\dot{x}_1 = u_1 & \qquad & \dot{y}_1 = v_1 & \qquad & \dot{z}_1 = w_1 \\
\dot{x}_2 = u_2 & \qquad & \dot{y}_2 = v_2 & \qquad & \dot{z}_2 = w_2 \\
\dot{x}_3 = x_2u_1 & \qquad & \dot{y}_3 = y_2 & \qquad & \dot{z}_3 = 1 + z_2w_1
\end{array}
\]
Equivalence maps: $(\bfx,\bfu) \leftrightarrow (\bfy,\bfv)$
\begin{eqnarray*}
\varphi( x,u,\dot{u},\ldots ) & = & (\ x_1x_2 - x_3,\ u_2,\ x_2,\ x_1u_2,\ \dot{u}_2,\ \ldots\ ) \\
\displaystyle\varphi^{-1}( y,v,\dot{v},\ldots ) & = & (\ v_1/y_2,\ y_3,\ y_3v_1/y_2 - y_1,\ \frac{y_2\dot{v}_1 - v_1v_2}{y_2^{\ 2}},\ y_2,\ \ldots\ )
\end{eqnarray*}
Verifying solutions:
\begin{eqnarray*}
\begin{array}{lcl}
	\begin{array}{rcl}
		\dot{y}_1 & = & \frac{d}{dt}(x_1x_2 - x_3) \\
		& = & x_1u_2 \\
		& = & v_1
	\end{array}
	& \qquad &
	\begin{array}{rcl}
		\dot{x}_1 & = & \D\frac{d}{dt}\left( \frac{v_1}{y_2} \right) \\
		& = & \D\frac{y_2\dot{v}_1 - v_1v_2}{{y_2}^2} \\
		& = & u_1
	\end{array}
	\\
	\begin{array}{rcl}
		\dot{y}_2 & = & \dot{u}_2 \\
		& = & v_2
	\end{array}
	& \qquad &
	\begin{array}{rcl}
		\dot{x}_2 & = & \dot{y}_3 \\
		& = & y_2 \\
		& = & v_2
	\end{array}
	\\
	\begin{array}{rcl}
		\dot{y}_3 & = & \dot{x}_2 \\
		& = & u_2 \\
		& = & y_2
	\end{array}
	& \qquad &
	\begin{array}{rcl}
		\dot{x}_3 & = & \D\frac{d}{dt}\left( \frac{y_3v_1}{y_2} - y_1 \right) \\
		& = & \D v_1 + y_3\frac{y_2\dot{v}_1 - v_1v_2}{{y_2}^2} - v_1 \\
		& = & x_2u_1
	\end{array}
\end{array}
\end{eqnarray*}

Equivalence map: $(\bfz,\bfw) \leftrightarrow (\bfy,\bfv)$
\begin{eqnarray*}
\psi( z,w,\dot{u},\ldots ) & = & (\ z_3 - z_1z_2,\ w_2,\ z_2,\ 1 - z_1w_2,\ \dot{w}_2,\ \ldots\ ) \\
\displaystyle\psi^{-1}( y,v,\dot{w},\ldots ) & = & (\ \frac{1-v_1}{y_2},\ y_3,\ y_1 + y_3\frac{1-v_1}{y_2},\ \frac{v_1v_2 - v_2 - y_2\dot{v}_1}{ y_2^{\ 2} },\ y_2,\ \ldots\ )
\end{eqnarray*}
Verifying solutions:
\begin{eqnarray*}
\begin{array}{lcl}
	\begin{array}{rcl}
		\dot{y}_1 & = & \D\frac{d}{dt}\left( z_3 - z_1z_2 \right) \\
		& = & 1 + z_2w_1 - w_1z_2 - z_1w_2 \\
		& = & v_1
	\end{array}
	& \, &
	\begin{array}{rcl}
		\dot{z}_1 & = & \D\frac{d}{dt}\frac{1-v_1}{y_2} \\
		& = & \D\frac{v_1v_2 - v_2 - y_2\dot{v}_1}{{y_2}^2} \\
		& = & w_1
	\end{array}
	\\
	\begin{array}{rcl}
		\dot{y}_2 & = & \dot{w}_2 \\
		& = & v_2
	\end{array}
	& &
	\begin{array}{rcl}
		\dot{z}_2 & = & \dot{y}_3 \\
		& = & y_2 \\
		& = & w_2
	\end{array}
	\\
	\begin{array}{rcl}
		\dot{y}_3 & = & \dot{z}_2 \\
		& = & w_2 \\
		& = & y_2
	\end{array}
	& &
	\begin{array}{rcl}
		\dot{z}_3 & = & \D\frac{d}{dt}\left( y_1 + y_3\frac{1-v_1}{y_2} \right) \\
		& = & \D 1 + y_3\frac{(v_1 - 1)v_2 - y_2\dot{v}_1}{{y_2}^2} \\
		& = & 1 + z_2w_1
	\end{array}
\end{array}
\end{eqnarray*}

Equivalence map: $(\bfx,\bfu) \leftrightarrow (\bfz,\bfw)$
\begin{eqnarray*}
\theta( x,u,\dot{u},\ldots ) & = & (\ \frac{1}{u_2} - x_1,\ x_2,\ \frac{x_2}{u_2} - x_3,\ -u_1 - \frac{\dot{u}_2}{u_2^{\ 2}},\ u_2,\ \ldots\ ) \\
\displaystyle\theta^{-1}( z,w,\dot{w},\ldots ) & = & (\ \frac{1}{w_2} - z_1,\ z_2,\ \frac{z_2}{w_2} - z_3,\ -w_1 -  \frac{\dot{w}_2}{w_2^{\ 2}},\ w_2,\ \ldots\ )
\end{eqnarray*}
Note that $\theta=\psi^{-1}\circ\phi$.

Verifying solutions:
\begin{eqnarray*}
\begin{array}{lcl}
	\begin{array}{rcl}
		\dot{x}_1 & = & \D \frac{d}{dt}\left( \frac{1}{w_2} - z_1 \right) \\
		& = & \D \frac{-\dot{w}_2}{{w_2}^2} - w_1 \\
		& = & u_1
	\end{array}
	& \qquad &
	\begin{array}{rcl}
		\dot{z}_1 & = & \D \frac{d}{dt}\left( \frac{1}{u_2} - x_1 \right) \\
		& = & \D \frac{-\dot{u}_2}{{u_2}^2} - u_1 \\
		& = & w_1
	\end{array}
	\\
	\begin{array}{rcl}
		\dot{x}_2 & = & \dot{z}_2 \\
		& = & w_2 \\
		& = & u_2
	\end{array}
	& \qquad &
	\begin{array}{rcl}
		\dot{z}_2 & = & \dot{x}_2 \\
		& = & u_2 \\
		& = & w_2
	\end{array}
	\\
	\begin{array}{rcl}
		\dot{x}_3 & = & \D \frac{d}{dt}\left( \frac{z_2}{w_2} - z_3 \right) \\
		& = & \D z_2 \left(-w_1 - \frac{\dot{w}_2}{{w_2}^2} \right) \\
		& = & x_2u_1
	\end{array}
	& \qquad &
	\begin{array}{rcl}
		\dot{z}_3 & = & \D \frac{d}{dt}\left( \frac{x_2}{u_2} - x_3 \right) \\
		& = & \D 1 + x_2 \left( -u_1 - \frac{\dot{u}_2}{{u_2}^2} \right) \\
		& = & 1 + z_2w_1
	\end{array}
\end{array}
\end{eqnarray*}

\end{proof}

\section{Conclusions}
\label{summary}

Wilkens showed that there are five equivalence classes of affine linear control systems with three state variables and two control variables under static equivalence. Below is a listing of how these classes combine using dynamic equivalence through one prolongation, using the static class representations presented in Elkin. Each equivalence class under dynamic equivalence is numbered. These nontrivial equivalences (or non-equivalences) are the work of this paper. \\

\begin{center}
\begin{tabular}{|c|c|c|c|}
\hline
$1$ & $\dot{x}_1 = u_1$ & $\dot{x}_1 = u_1$ & $\dot{x}_1 = u_1$ \\
 & $\dot{x}_2 = u_2$ & $\dot{x}_2 = u_2$ & $\dot{x}_2 = u_2$ \\
 & $\dot{x}_3 = x_2$ & $\dot{x}_3 = x_2u_1$ & $\dot{x}_3 = 1 + x_2u_1$ \\
\hline
$2$ & $\dot{x}_1 = u_1$ & & \\
 & $\dot{x}_2 = u_2$ & & \\
 & $\dot{x}_3 = 0$ & & \\
\hline
$3$ & $\dot{x}_1 = u_1$ & & \\
 & $\dot{x}_2 = u_2$ & & \\
 & $\dot{x}_3 = 1$ & & \\
\hline
\end{tabular}
\end{center}

Future avenues of research into the classification of affine linear control systems under dynamic equivalence include looking at higher order equivalences ($J$ and/or $K>0$) as well as increasing the number of state and control variables. One obstacle to overcome with higher order equivalences and larger numbers of variables is that, unlike the case presented here where a unique $\bfS$ exists, the problem quickly splits into many cases with different $\bfS$. In addition, this method relies on the fact that affine linear systems in this dimension have already been classified under static equivalence, and the static equivalence problem for affine control systems has only been completed in a few low-dimensional cases. Nevertheless, the further exploration of this decomposition may still yield new insights into the phenomenon of dynamic equivalence in general.

\section{Acknowledgements}

I would like to thank my thesis adviser, Jeanne Clelland, for all her guidance, both research related and not. I would also like to thank George Wilkens who introduced me to the topic of dynamic equivalence and helped me along the way.

\end{document}